\documentclass[twoside]{article}
\usepackage{eurosym}
\usepackage[utf8]{inputenc}
\usepackage{graphicx}
\usepackage{latexsym}
\usepackage{amsmath}
\usepackage{amsfonts}
\usepackage{amssymb}
\usepackage{hyperref}
\usepackage[usenames]{color}

\newcommand{\disp}{\displaystyle}
\newcommand{\proof}{\par
	\noindent {\sc Proof}.\quad}

\newcommand{\rA}{{\rm A}}

\newcommand{\rB}{{\rm B}}
\newcommand{\bB}{{\bf B}}

\newcommand{\cC}{{\cal C}}
\newcommand{\rC}{{\rm C}}

\newcommand{\rD}{{\rm D}}

\newcommand{\rE}{{\rm E}}

\newcommand{\rF}{{\rm F}}

\newcommand{\rG}{{\rm G}}
\newcommand{\cG}{{\cal G}}

\newcommand{\rH}{{\rm H}}

\newcommand{\rI}{{\rm I}}

\newcommand{\rK}{{\rm K}}

\newcommand{\rL}{{\rm L}}

\newcommand{\cM}{{\cal M}}
\newcommand{\rM}{{\rm M}}
\newcommand{\rN}{{\rm N}}
\newcommand{\NN}{{\mathbb{N}}}

\newcommand{\bn}{{\bf n}}
\newcommand{\rn}{{\rm n}}

\newcommand{\rP}{{\rm P}}

\newcommand{\RR}{{\mathbb{R}}}

\newcommand{\rR}{{\rm R}}

\newcommand{\rS}{{\rm S}}
\newcommand{\rT}{{\rm T}}

\newcommand{\rU}{{\rm U}}
\newcommand{\rV}{{\rm V}}

\newcommand{\rW}{{\rm W}}

\newcommand{\rX}{{\rm X}}

\newcommand{\rY}{{\rm Y}}
\newcommand{\ry}{{\rm y}}

\newcommand{\fin}{\hfill\mbox{$\quad{}_{\Box}$}}
\newcommand{\fineq}{\vspace{-.75cm$\fin$}\par\bigskip}
\newcommand{\fineqnum}{\vspace{-.4cm$\fin$}\par\bigskip}

\newcommand{\pe}[2]{\disp{\langle #1,#2\rangle}}

\newtheorem{theo}{\bf \sffamily Theorem}
\newtheorem{prop}{\bf \sffamily Proposition}
\newtheorem{rem}{\bf \sffamily Remark}
\newtheorem{lemma}{\bf \sffamily Lemma}
\newtheorem{definition}{\bf \sffamily Definition}
\newtheorem{exam}{\bf \sffamily Example}
\begin{document}
	
\title{{\Large \bfseries\sffamily Controlled boundary explosions: dynamics after blow-up for some semilinear problems with global controls}
\author{\bfseries\sffamily A.C. Casal, G. D\'{\i}az, J.I. D\'{\i}az, J.M. Vegas }
			\date{\qquad \boxed{\bf Revised.~April, 26}}
\date{}
}	
\maketitle
\begin{center}
	\textit{Dedicated to Juan Luis Vázquez on occasion of his 75th birthday}
\end{center}
\begin{abstract}
The main goal of this paper is to show that the blow up phenomenon (the explosion of the $
\rL^{\infty }$-norm) of the solutions of several classes of evolution problems
can be controlled by means of suitable global controls $\alpha (t)$ ($i.e.$ only
dependent on time) in such a way that the corresponding solution be well
defined (as element of $\rL_{loc}^{1}(0,+\infty :\rX)$, for some functional
space $\rX$) after the explosion time. We start by considering the case of an
ordinary differential equation with a superlinear term and show that the controlled 
explosion property holds by using a delayed control (built through the solution of the problem and by generalizing the {\em nonlinear variation of constants formula}, due to V.M.
Alekseev in 1961, to the case of {\em neutral delayed equations} (since
the control is only in the space $\rW_{loc}^{-1,q\prime }(0,+\infty :\RR
)$, for some $q>1$)$.$ We apply those arguments to the case of an evolution semilinear problem in which the differential equation is a semilinear elliptic equation with a
superlinear absorption and the boundary condition is dynamic and involves
the forcing superlinear term giving rise to the blow up phenomenon. We prove that,
under a suitable balance between the forcing and the absorption terms, the
blow up takes place only on the boundary of the spatial domain which 
here is assumed to be a ball $\rB_{\rR}$ and for a constant as initial datum.

\end{abstract}

\section{Introduction}

It is well known than one of the more relevant qualitative behaviors of
nonlinear evolution problems is the possibility to get the finite time blow-up
of the $\rL^{\infty}$-norm of the solution of suitable parabolic problems. Without any aims to be
exhaustive, we mention as general references the books \cite{ Galaktionov-Vazquez Libro, Hu, Quitner-Souplet} as well as the revision made in  \cite{Brezis-Cacenave}.

The main problem we will consider in this paper is a semilinear dynamic boundary condition for an elliptic diffusion absorption equation on a ball $\bB_{\rR}\subset \RR^{\rN},~\rN>1$,
\begin{equation}
\rP(\alpha)\quad \left\{
\begin{array}[c]{ll}
-\Delta u+\cG(u,\alpha )=0 & \quad\hbox{in $\bB_{\rR}\times (0, +\infty  )$},\\[.15cm]
\dfrac{\partial u }{\partial t}+\dfrac{\partial u }{\partial \bn}=\rF(u,\alpha ) & \quad \text{on }\partial \bB_{\rR}\times (0,\infty), \\[.2cm]
u(x,0)\equiv u_{0}, & \quad x\in\partial\bB_{\rR},
\end{array}
\right.  \label{eq:controlledproblem}
\end{equation}
where $\rF,\cG:\RR\times \RR \rightarrow \RR$ are two locally Lipschitz real functions, $u_{0}$ is a positive constant and the control $\alpha $ is global, in the
sense that $\alpha $ is a {\em purely time-depending} scalar function $
\alpha (t)$. The crucial fact in our study is the assumption that, in the absence of control $\alpha (t)\equiv 0$, functions $\rF(u,0)$ and $\cG(u,0)$ are superlinear. More exactly, we will assume that 
\begin{equation}
\rF(u,0)=\lambda f(u),  \label{Hypo alfa cero}
\end{equation}
where $\lambda >0$ is a given parameter and $f:\RR\rightarrow 
\mathbb{R}$ \ is a locally Lipschitz real function "superlinear near 
infinity" in the sense that 
\begin{equation}
\int_{r_{0}}^{\infty }\frac{ds}{f(s)}<\infty,   \label{Hypo superlinear}
\end{equation}
\noindent for some $r_{0}>0$,  and that 
\begin{equation}
\mathcal{G}(u,0) =g(u),  \label{Hypo alfa cero g}
\end{equation}
where $g:\RR\rightarrow 
\mathbb{R}$ \ is a locally Lipschitz real function "superlinear near the infinity" in the Keller--Osserman sense:  
\begin{equation}
\int_{r_{0}}^{\infty}\dfrac{ds}{\sqrt{\rG(s)}}<\infty, 
\label{eq:KellerOssermanG}
\end{equation}	
for some $r_{0}>0$, where $\displaystyle \rG(s)=\int^{s}_{0}g(s)ds$. Some examples intensively studied in the literature correspond to the cases $f(u)=e^{u},$ $f(u)=(k+u)^{p}$ for some $p>1$ and $k\geq 0,$ and similar choices also for $g$.
We will need also the assumption that the forcing term dominates over the absorption one, in the sense that 
\begin{equation}
	\liminf_{\tau\rightarrow\infty}\dfrac{f(\tau)}{\sqrt{2\rG\big (\tau\big )}}\in (0,+\infty ].
		\label{eq:dominationbalance}
\end{equation}
We will see later that this domination assumption has different consequences according to the above limit is infinite (strong domination) or positive (weak domination). In any case,  we note that, under this domination assumption, the Keller--Osserman condition implies the superlinear condition \eqref{Hypo superlinear} for function $g$.

We will denote by $u^{\alpha}$ the solution of $\rP(\alpha)$, and thus $u^{0}$ will represent the solution of $\rP(0)$ (the problem without any control). We will show that solutions $u^{\alpha}$ attain their first blow-up only on the boundary of the spatial domain $\partial \bB_{\rR}$. This must be done in contrast with what happens in the usual case of semilinear parabolic problem with Dirichlet boundary conditions (for which the first blow-up may take place in a single point), or Neumman or Robin boundary conditions (with  the first blow-up takes place taking place in the whole domain $\Omega$). As far as we know, problem $\rP(\alpha)$ was not considered in the previous literature and will be the object of a systematic study in Section~\ref{section:completerecuperation}. After that, the second goal of this paper is to show that, given any small $\varepsilon>0$ , it is possible to find a suitable control  $\alpha (t)$ which allows to ensure that (after, perhaps, a small modification of the superlinear terms) the corresponding solution can be continued beyond the finite time blow-up of its $\rL^{\infty}$-norm (which is the same than the one of the solution corresponding to no control case $\alpha (t)=0,~\rF(u,0)=\lambda f(u)$ and $\mathcal{G}(u,0) =g(u)$). 
This property should be contrasted with the {\em complete blow-up phenomenon}, which holds in  most of the usual superlinear problems (see, $e.g.$, \cite{Baras-Cohen} and \cite{Galaktionov-Vazquez Libro}). 
As we will see, the key tool in our study will be the comparison principle
jointly with the previous study of a simpler case: {\em the superlinear ordinary differential equation}
$$
\rP_{\rF(\cdot,\alpha )} \quad \left\{ 
\begin{array}{lc}
\dfrac{du^{\alpha}}{dt}(t)=\rF(u^{\alpha},\alpha) & \text{in } (0,+\infty  ), \\ [.15cm]
u^{\alpha}(0)=u_{0}, & 
\end{array}%
\right.   
$$
with $\rF(u,\alpha)$ satisfying \eqref{Hypo alfa cero} and the structure condition given below (see in \eqref{hyp quasibangnbangF}).

The property we will study in this paper can be stated in the following terms:

\begin{definition}
We say that the solution $u^0(\cdot,t)$ of problem $\rP(0)$ {\rm (}respectively $\rP_{\rF(\cdot,0)}${\rm )}, with no control $\alpha=0$ {\rm (} with blow-up time $\rT_{\infty}(u^0)${\rm )}, has a "controllable explosion" if 
\par
\noindent i) for any given $\varepsilon>0$ we can find a
continuous deformation and an extension of the trajectory
$u^0(t,.)$ on the interval $[0,\rT_{\infty}(u^0)-\varepsilon)$, by $u^{\alpha}(\cdot,t)$, solution of the associate perturbed control problem defined by replacing $\lambda f(u)$ by $\rF(u,\alpha)$ and $g$ by $\cG(u,\alpha )$,
for a suitable control $\alpha$, such that 
$u^{\alpha}(\cdot,t)$ also blows-up at the same time $\rT_{\infty}(u^{\alpha})=\rT_{\infty}(u^0)$,
\par
\noindent ii) $u^{\alpha}(\cdot,t)$ can be extended beyond $\rT_{\infty}(u^{\alpha})$ as a solution which is in the space
$\rL_{loc}^{1}(0,+\infty :\rX)$ {\rm (}i.e. in $\rL^{1}(0,\rT:\rX)$, for any $\rT\in (0,+\infty )${\rm )} where $\rX=\rL^{\infty}(\bB_{\rR})$ {\rm (}respectively $\rX=\RR ${\rm )}.
\end{definition}

As a matter of fact, we will prove a stronger conclusion for  the corresponding solution $u^{\alpha }(t)$ of the problems $\rP(\alpha)$ and $\rP_{\rF(.,\alpha )}$: we will prove that, in fact, $u^{\alpha
}(t)=u^{0}(t)$ for any $t\in [0,\rT_{\infty}(u^{0})-\varepsilon ]$. Notice that we use the notation $ \rT_{\infty}(u^{0})$ since the blow up time depends not only on the initial datum $u_{0}$ but also on other parameters and data of the problem. In the rest of paper we will simplify the notation by writing  $\rT_{\infty}(u^{0})=\rT_{\infty}$.

This philosophy of controlling in order to have the recuperation after the explosion was initiated in the previous papers by the authors dealing with some special multiplicative control of delayed feed-back form for some ordinary and partial differential equations (\cite{CDV2008, CaDiVe2013} and
\cite{CaDiVe}). 

This approach applies even to certain linear delayed problems. Here we will extend, and improve, this approach by considering nonlinear terms leading to the explosion on the boundary of the spatial domain. We will make mention also of some other problems at the end of the paper (see Remark \ref{Other problems}). We will prove that the recuperation after the blow up time arises for controls $\alpha (t)$ with a quasi-bang-bang structure
\begin{equation}
\rF(u,\alpha )=\rS\lambda f_{\rM_{\epsilon}}(u)+\alpha,
\label{hyp quasibangnbangF}
\end{equation}
$$
f_{\rM_{\epsilon}}(u)=\left\{ 
\begin{array}{ll}
f(u) & \text{if }0\leq u\leq \rM_{\epsilon}, \\ [.1cm]
f(\rM_{\varepsilon}) & \text{if }u>\rM_{\epsilon},%
\end{array}%
\right. 
$$  
for some $\rM_{\varepsilon}>0,$ and with $\rS\in \hbox{\rm sign}^{\pm }(\alpha )$, where 
$$
\hbox{\rm sign}^{\pm }(\alpha )=\left\{ 
\begin{array}{cc}
1 & \text{if }\alpha \geq 0, \\ 
-1 & \text{if }\alpha \leq 0.%
\end{array}%
\right. 
$$
Notice that $\hbox{\rm sign}^{\pm }(0)=\{-1,+1\}.$ Moreover we will see that $\alpha(t)= 0$ if $t\in (0,\rT_{\infty}-\varepsilon)$ and $\alpha(t)> 0$ in $t\in \big (\rT_{\infty }-\varepsilon,\rT_{\infty }\big )$. We will show that in some cases, the "effective control" becomes $\rS(\alpha)$ and it is purely of bang-bang type and no truncation is needed
\begin{equation}
\rF(u,\alpha)=\rS(\alpha)\lambda f(u),
\label{hyp purebangnbang}
\end{equation}
with $\rS(\alpha)\in \hbox{\rm sign}^{\pm }(\alpha )$. This is the case when, for instance, $f(u)=(k+u)^{p}$ with $p>2$ and $k\geq 0$ (see Remark \ref{rem:bang-bang}).
The dependence on $\alpha$ of the function $\cG(u,\alpha )$ is much weaker since we will need only to truncate function $g$ if $\alpha\neq 0$: 
\begin{equation}
\cG(u,\alpha )=g_{\rM_{\epsilon}}(u)  \quad \text{if } \alpha\neq 0,\\
\label{hyp quasibangnbangGNew}
\end{equation}
where
$$
g_{\rM_{\epsilon}}(u)=
\left\{ 
\begin{array}{ll}
g(u) & \text{if }0\leq u\leq \rM_{\varepsilon}, \\ [.15cm]
g(\rM_{\varepsilon}) & \text{if }u>\rM_{\varepsilon}, 
\end{array} 
\right .  
$$
for some $\rM_{\varepsilon}>0$.

It is clear that the recuperation of the solution after the blow up time $
\rT_{\infty}$ requires that $u_{t}^{\alpha}(t)<0$ for $t\in (\rT_{\infty },\rT_{\infty
}+\delta)$, for some $\delta >0$ and thus it is natural to assume the quasi-bang-bang structure and to take  $\alpha (t)<0$ for $t\in (\rT_{\infty },\rT_{\infty }+\delta)$. The
quasi-bang-bang structure of the controls implies that the initial superlinear
forcing term $\lambda f(u)$, on the interval $[0,\rT_{\infty })$ becomes later a 
{\em superlinear absorption} term, $-\lambda f(u)$, at least in a short period after 
 $\rT_{\infty }$ and the problem has infinity as initial value. The possibility to solve nonlinear parabolic problems with an infinite initial value, in presence of a {\em 
superlinear absorption term } was already proved for some special parabolic
problems (see, $e.g.$, \cite{BaDiDi} and its references). Here we must take into account the possibility of truncating $f(u)$ and the presence of a negative control $\alpha(t)$.

From the point of view of Control Theory, one of the pioneering works on
control for blow-up problems for nonlinear parabolic equations with a forcing 
term was the book by J.L. Lions \cite{Lions} (see also \cite{Diaz-Lions, Diaz-Lions 2, Coron-Trelat, Amman-Quittner, Balch, FdezCara-Zuazua, Porreta-Zuazua} and the references therein). In these, and many other works, the goal was to avoid the occurrence of the blow-up phenomenon by means of
suitable controls (the case of controls given by measures was considered in \cite{Amman-Quittner}). The possibility to choose the blow-point time and points
were considered in \cite{Merle} and
\cite{Cacenave-Matel-Zhao} for the nonlinear wave equation. The approximate controllability for the case of dynamic boundary conditions leading to global
solutions was considered in \cite{BeDiVr}. As far as we know, no previous attempt to control problem $\rP(\alpha)$, searching a continuation dynamics after the same blow-up time than the one of the solution without control, was
considered before. Notice the structure of our control problem is not entirely conventional since we allow to consider the case in which f is sign reversed (and truncated). Moreover, our control will be built as a solution
of an auxiliary, singular ODE with delay. This point of view is in contrast with many of the above mentioned control papers on control of semilinear partial differential parabolic problems in which the main non-linearity is kept, the control is additive and many times with a localized spatial support of the control for a given time horizon $\rT$ (instead, $\rT=\infty$), etc. 
A quite complete list of references dealing with nonlinear problems with dynamic boundary conditions, starting already in 1901, can be found, $e.g.$, in
the survey papers \cite{BeDiVr} and \cite{Bandle-Reichel}. The study of the
special case in which only the nonlinear dynamic boundary conditions is the
origin of blow-up phenomena was considered in \cite{Kirane} and later by
several other authors (see, $e.g.$, \cite{KiraneNabanaPokhozhaev}) but for different elliptic equations on the spatial domain. Notice that this is a different situation to the case in
which there is a nonlinear parabolic equation with a source term jointly with
a dynamic boundary condition (see, $e.g.$,\cite{Amman-Fila, Bandle-Reichel, Vazquez-Vitillaro}). In all the cases, the blow-up takes place also on the boundary, as it is the case of a
nonlinear parabolic equation with a source term jointly with a static possibly
nonlinear Robin type boundary condition (see, $e.g.$, \cite{Levine-Payne,Lopez-Wolansky} and the survey \cite{Fila-Filo}).

We point out that in fact the blow-up phenomenon only occurs for large enough
initial data (or large values of the parameter $\lambda$) since otherwise the solution is well defined for any $t>0$ and converges (as $t$ goes to infinity) to a solution of the stationary problem. This is well known for the usual semilinear parabolic problem. For problem $\rP(\alpha)$ it can be proved as an easy modification of some previous papers in the literature (see, $e.g.$
\cite{Fila-Quitner} and \cite{Arrieta}). This is the reason why we shall always assume an additional condition 
\begin{equation}
\hbox{$u_{0}$ or/and $\lambda $  are large enough}.
\label{Dato inicial grande}
\end{equation}
The key idea in this paper is to start, in Section 2, by proving that problem $\rP_{\rF(\cdot,\alpha )}$ admits the controlled explosion property and then to use similar ideas, in Section 3, for the case of problem $\rP(\alpha)$. 

For the case of $\rP_{\rF(\cdot,\alpha )}$, we will construct, in Section 2, the suitable control $\alpha (t)$
as a changing sign delayed term of the form $\rB^{\prime }(t)\ry(t-\tau )$, for a suitable function $\rB(t)$, where $\ry(t)$ is the solution of the
auxiliary "neutral delayed ordinary differential equation"
\begin{equation}
\left\{ 
\begin{array}{l}
\dfrac{d}{dt}\left[ \ry(t)-\rB(t)\ry(t-\tau )\right] =\lambda f_{\rM_{\epsilon}}(y)-\rB(t)\dfrac{d}{
dt}\left[ \ry(t-\tau )\right] ,\quad t>0 \\[0.2cm]
\ry(\theta )=u^{0}(\theta ),\quad 0\leq \theta \leq \rT_{\infty }-\varepsilon ,
\end{array}%
\right.   \label{Neutral ODE}
\end{equation}%
with the history initial condition
$$
\ry^{0}(\theta )=u^{0}(\theta )\text{ for any }\theta \in [ 0,\rT
_{\infty }-\varepsilon ],
$$
where $u^0(t)$ is the solution of problem $\rP_{\rF(\cdot,0)}$ with no control ($i.e.$ $\alpha=0$). 
\par
We emphasize that the good control will change sign with time and that
$$
\alpha (t)=\rB^{\prime }(t)y(t-\tau ).
$$
Moreover, as we will see, $\alpha \notin \rL_{loc}^{1}(0,\rT_{\infty} :\RR)$
but $\alpha \in \rW_{loc}^{-1,q\prime }(0,\rT_{\infty} :\RR)$, the dual
space of $\rW_{0,loc}^{1,q}(0,\rT_{\infty} :\RR)$, for some $ q>1.$
\par
In fact, the previous paper \cite{CDV2008} was devoted to a class of delayed
problems (as for instance problem (\ref{Neutral ODE}) with $\lambda =0$) and so the searched control in
that case was the function  $\rB^{\prime }(t)$, so the bang-bang term, $\hbox{\rm 
sign}^{+}(\alpha (t))$, was not needed. 
Here we will apply the main philosophy of the results of \cite{CaDiVe} (in which a refined {\em nonlinear
variation of constants formula}, initially due to Alekseev \cite{Alekseev},
was a crucial tool), to prove that $\ry\in \rL^{1}\big (0,\rT_{\infty }\big ),$ and that the extended solution satisfies problem $\rP_{\rF(\cdot,\alpha )}$ thanks to the quasi-bang-bang term (if $\lambda>0$) and an argument of time reflection over $\rT_{\infty }$ for the consideration of the corresponding superlinear absorption equation with infinity as initial datum. Then we will continue the solution to the whole interval $(0,\infty)$, by periodicity. This will be detailed in Section 2. 


We start Section \ref{section:completerecuperation} by developing the study of the problem $\rP(0)$ , $i.e.$, without any control.  We shall prove (Theorem \ref{theo:blowupsuperdomination})
that if the elliptic equation is of superlinear type ($i.e.,~g(u)=u^{m}$ with
$m>1$), as well as the forcing  term on the boundary (assumed, for instance, as $~f(u)=u^{p}$ with $p>1$), then, if the absorption rate is lower than the forcing one (this is the assumption (\ref{eq:dominationbalance}); so that $p\ge (m+1)/2>1$), in absence
of any control ($\alpha\equiv0$), the corresponding solution of  \eqref{eq:controlledproblem} has a finite blow up time $\rT_{\infty}$, so that 
$$
\left\{
\begin{array}[c]{ll}
-\Delta u^{0}\big (\cdot,\rT_{\infty}\big )+g\big (u^{0}\big (\cdot,\rT_{\infty}\big )\big )=0 & \quad\hbox{in $\bB_{\rR }$},\\[0.175cm]
u ^{0}\big (\cdot,\rT_{\infty}\big )=+\infty & \quad\hbox{on $\partial \bB_{\rR}$.}
\end{array}
\right. 
$$
Thus, at the explosion time $u^0$ coincides with the unique large solution, $\rU_{\infty}^{\bB_{\rR}}$, of the associate elliptic problem (see \eqref{eq:largesolution} below). Since we have the inequality 
$$
u^{0}(x,t) \le \rU_{\infty}^{\bB_{\rR}}(x)<+\infty,\quad x\in\bB_{\rR},~0<t,
$$
then the explosion is only possible on $\partial \bB_{\rR}$ after $\rT_{\infty}$. It make sense thanks to the domination assumption (\ref{eq:dominationbalance}) which in the case of powers it corresponds to the condition $p>\dfrac{m+1}{2}>1$. We will obtain also several time estimates on the behaviour of solutions near the finite blow up time $\rT_{\infty}$. We will pay attention also to the limiting (weak domination) case $p=(m+1)/2$ (see Theorem \ref{theo:invariance}) of the domination assumption. In order to better illustrate the behaviour near the finite blow up time $\rT_{\infty}$, we consider in this case  a self-similar solution corresponding to the spatial domain given by the hyperplane $\RR^{\rN-1}\times \RR_{+}$. In the case of general nonlinear terms, $f$ and $g$, some additional technical assumptions are required, as we will indicate below. For the controlled problem $\rP(\alpha)$ we will show, again, that it is possible to choice a control $\alpha(t)$, now acting on the boundary of the spatial domain $\partial\rB_{\rR}$, such that the corresponding solution $u ^{\alpha}(t)$ satisfies the properties indicated in Definition 1 (see Theorem \ref{Thm PIII}).

\section{The control for the complete recuperation after the blow up time for problem $\rP_{\rF(\cdot,\alpha )}$}
\label{section:completerecupartionalpha}

It is well known that the simpler and illustrative example of dynamical system for which there is blow up of solutions is the problem  $\rP_{\rF(\cdot,\alpha )}$ (before the blow up time)
$$
\left \{
\begin{array}{l}
	\dfrac{d u^{0}}{dt} (t)=\lambda f\big (u^{0}(t)\big),\quad t>0,\\[.2cm]
	u^{0}(0)=u_{0}>0.
\end{array}
\right .
$$
where $\lambda$ is a positive constant and $f:~\RR_{+}\rightarrow \RR_{+}$ is a continuous function. Since 
$$
\dfrac{\dfrac{d u^{0}}{dt}(t)}{f\big (u^{0}(t)\big)}=\lambda ,
$$ 
assuming condition \eqref{Hypo superlinear} we deduce
$$
\int^{u^{0}(\widehat{t})}_{u^{0}(t)}\dfrac{d s}{f(s)}=\lambda \big (\widehat{t}-t\big),\quad 0<t<\widehat{t}.
$$
We note that  \eqref{Hypo superlinear} enables us to consider the decreasing function
$$
\Phi(r)=\int^{+\infty}_{r}\dfrac{ds}{f(s)},\quad s>0,
$$
with $\Phi(+\infty)=0$. This function will be used systematically in Section \ref{section:completerecuperation}. We have that
$$
t\mapsto \Phi\big (u^{0}(t)\big )+\lambda t
$$
is a constant function.  In particular, we may define $\rT=\dfrac{\Phi (u_{0})}{\lambda}$ for which
$$
\left \{
\begin{array}{ll}
	u^{0}(t)=\Phi ^{-1}\big (\lambda (\rT-t)\big )<+\infty,&\quad 0<t<\rT,\\ [.2cm]
	\disp u^{0}(\rT^{-})\doteq \lim_{t\nearrow \rT}u^{0}(t)=+\infty.
\end{array}
\right .
$$
Notice that with the notation of the Introduction $\rT_{\infty}=\dfrac{\Phi (u_{0})}{\lambda}$ in this problem and that here $\lambda >0$ is arbitrary.
\begin{rem}\rm For the power case $f_{p}(s)=s^{p},~p>0,$ the condition \eqref{Hypo superlinear} corresponds to $p>1$. Then $\Phi_{p}(s)=\dfrac{1}{p-1}\dfrac{1}{s^{p-1}}$ and 
$$
u^{0}(t)=\dfrac{1}{(p-1)^{\frac{1}{p-1}}}\dfrac{1}{\big (\lambda (\rT_{\infty}-t)\big)^{\frac{1}{p-1}}}, \quad 0<t<\rT_{\infty}\doteq \dfrac{1}{\lambda (p-1)u_{0}^{p-1}}.
$$
\fineq
\label{eq:powerlikef}
\end{rem}
\begin{rem}\rm It is clear that if 
$$
\widehat{f}(s)\ge f(s)\quad \hbox{for large $s$}
$$
then one deduces that if $f$  verifies  \eqref{Hypo superlinear} the same happens with $\widehat{f}$. In particular, any function $\widehat{f}(s)\ge s q(s)$, for large $s$,  verifying 
$$
\liminf_{s \rightarrow \infty} \dfrac{q(s)}{s^{\gamma}}\in (0,+\infty]\quad 
\hbox{for some $\gamma >0$,}
$$
satisfies \eqref{Hypo superlinear}. For instance, we may choose $q(s)\ge \big (\log s\big )^{\gamma},~\gamma\ge 1$, or $q(s)\ge \log (\log(\cdots \log (s))).\fin$
\label{rem:exemplesf}
\end{rem}
\begin{rem}\rm From Remark \ref{rem:exemplesf} it follows that if we assume the property
\begin{equation}
\dfrac{f(s)}{s^{\alpha }}\quad \hbox{is increasing for large $s$},
\label{eq:monotocitysuperlinearf}
\end{equation}
for some $\alpha>1$, then the assumption \eqref{Hypo superlinear} is satisfied. Moreover, if $\nu>1$, we deduce that
$$
\Phi(r)=\int^{+\infty}_{r}\dfrac{ds}{f(s)}=\nu \int^{+\infty}_{\frac{r}{\nu}}\dfrac{d\widehat{s}}{f(\nu\widehat{s})}\le
\nu^{1-\alpha}\Phi\big (\nu^{-1}r\big )\quad \hbox{for large $r$}.
$$
Therefore, the change of variables $\zeta=\Phi (r)$ implies that
\begin{equation}
\nu\Phi^{-1} \left (\nu^{\alpha -1}\zeta\right )\ge \Phi^{-1}(\zeta),\quad \hbox{for small $\zeta$},
\label{eq:technicallityf}
\end{equation}
for $\nu>1$. A similar argument enables us to obtain
\begin{equation}
\nu\Phi^{-1} \left (\nu^{\alpha -1}\zeta\right )\le \Phi^{-1}(\zeta),\quad \hbox{for small $\zeta$},
\label{eq:technicallityf2}
\end{equation}
for $\nu<1.\fin$ 
\label{rem:technicalityf}
\end{rem}
\par

The main result of this section devoted to the nonlinear ordinary differential equation is the following:
\begin{theo} \label{teo:theo1} Assume  $f$ locally Lipschitz continuous and superlinear. Then, for any $u_{0}>0$ the blowing up trajectory $u^{0}(t)$ 
of the associated problem $\rP_{\rF(\cdot,0)}$ has a controlled
explosion {\rm (}in the sense of Definition 1{\rm )} by means of the control problem $\rP_{\rF(.,\alpha )}$ for a suitable $\alpha \in \rW_{loc}^{-1,q\prime }(0,\rT_{\infty} :\RR)$, for some $ q>1.$
Moreover, given $\varepsilon>0$, if $u^{\alpha }\in\rL_{loc}^{1}(0,+\infty)$ is the solution of $\rP_{\rF(\cdot,\alpha )}$ corresponding to the built control $\alpha (t)$, then $u^{\alpha }(t)$ coincides with the solution with no control  $u^{\alpha
}(t)=u^{0}(t)$  for any  $t\in [0,\rT_{\infty}-\varepsilon ]$ and $u^{\alpha }(t)$  also blows-up at the time $\rT_{\infty}$ corresponding to $u^{0}(t)$. In
addition, if $\alpha (t)$ is the searched
control {\rm (}in the sense of Definition 1{\rm )} then the corresponding solution $
u^{\alpha }(t)$ satisfies that $u^{\alpha }(t)>0$ for any $t>0.$
\end{theo}

Our main tools, and the strategy, of the proof are the following: we start by taking a delayed feedback control (in the spirit of \cite{CDV2008}), on the interval $[0,\rT_{\infty}]$, with $\rT_{\infty}=\rT_{\infty}(u^{\alpha})$,  in order to the associated solution to be in 
$\rL^{1}(0,\rT_{\infty})$. To this end we will apply the so called \textit{nonlinear variation of constants formula} to the problem with the truncated function $f_{\rM_{\varepsilon}}(u)$. After that, we pass to the consideration of the corresponding superlinear absorption problem with infinity as initial datum, for $t\in [\rT_{\infty},2\rT_{\infty})$. Finally, a periodicity argument allow to extend the solution to the interval $t\in (2\rT_{\infty},+\infty)$ (see Figure 1 below).
\par
Concerning the \textit{nonlinear variation of constants formula} we recall that it was first established in the
literature for\textit{\ }nonlinear terms $h$ of class $\cC^{2}$ (see Alekseev \cite
{Alekseev}, Laksmikantham and Leela \cite{LaksLeela}). Here we will prove that the formula holds also for Lipschitz functions $h$ \ (which at this stage can be assumed to be in fact globally
Lipschitz) and with a very general perturbation term (which in fact can be a
multivalued term). Given a family of maximal
monotone operators $\beta (t,y)$, on the space $\rH=\RR^{d}$, with $
\beta (\cdot,t)\in \rL_{loc}^{1}(0,+\infty :\RR^{d}),$
we consider the perturbed problem
$$
\rP^{\ast }(h,\beta ,\xi )=\left\{ 
\begin{array}{l}
\dfrac{d\ry}{dt}(t)+\beta (\ry(t),t)\ni  \ h(\ry(t)),\text{ in }\RR^{d}, \\ [.15cm]
\ry(t_{0})=\xi .
\label{eq:multivalued}
\end{array}
\right.
$$
We know that once that $h$ is globally Lipschitz function, the solutions of $
\rP(h,\beta ,\xi )$ are well defined, as absolutely continuous functions on $
[0,\rT],$ for any given $\rT>0$ (this is an easy consequence of the general
theory: see \cite{BreOMM} for the autonomous case, and \cite{Ya}, and its references, for the generalizations to the case of $\beta $ depending on $t$). 

Now, we reformulate the trajectory $\ry^{0}(t)$ of \eqref{eq:multivalued} with $\beta \equiv 0$ in more general terms (by modifying the initial time and the initial
condition). So, we define  $\ry^{0}(t)=\phi (t,t_{0},\xi )$, with $\phi (t,t_{0},\xi )$ the
unique solution of the ODE
$$
\rP^{\ast }(h,0,\xi )=\left\{ 
\begin{array}{l}
\ry^{\prime }(t)= h(\ry(t))\text{ \ in }\RR^{d}, \\ [.15cm]
\ry(t_{0})=\xi .%
\end{array}%
\right.
$$
We introduce the formal notation $\Phi (t,t_{0},\xi )=\partial _{\xi }\phi
(t,t_{0},\xi ),$ where $\partial _{\xi }$ denotes the partial
differentiation. Then we shall prove:

\begin{theo} \label{teo:theo2} The flow map $\phi $ is 
Lipschitz continuous, $\Phi $ is absolutely continuous
and the solution $\ry(t)$ of the "perturbed"
problem $\rP^{\ast }(h,\beta ,\xi )$ has the integral representation
$$
\ry(t)=\ry^{0}(t)-\int_{t_{0}}^{t}\Phi (t,s,\ry(s))\beta (s,\ry(s))ds\quad \hbox{for any }t\in [0,\rT],
$$
where $\ry^{0}(t)=\phi (t,t_{0},\xi )$ is the solution of
the unperturbed problem $\rP^{\ast}(h,0,\xi )$.
\end{theo}
In the above formula we used, for simplicity, the notation corresponding to the case in which  $\beta (\cdot,t)$ is single-valued, but a suitable similar expression
can be formulated if $\beta (\cdot,t)$ is multivalued. As a matter of fact, we will also generalize (for the case $d=1$) the Alekseev's formula to the case in which the perturbation $\beta (t,y(t))$ of the equation is an element in the space $\rW^{-1,q\prime }(0,\rT:\RR)$.

\subsection{Proof of Theorem 1 assuming Theorem 2}

We assume, for a while, that Theorem 2 holds. The proof of Theorem 1 can be divided in different steps.
\par
\noindent {\em Step 1}: $t\in[0,\rT_{\infty}]$. Let us develop the indicated  strategy on the initial interval $[0,\rT_{\infty}]$. Given $\varepsilon >0$, we define $
\tau =\rT_{\infty}-\varepsilon $ and $\rM_{\varepsilon }=u^{0}(\rT_{\infty}-\varepsilon )$ (this explains the dependence on $\varepsilon$ of the truncation parameter $\rM_{\varepsilon }$. Notice that then $f_{\rM_{\varepsilon }}(u^{\alpha}(t))=f(u^{\alpha}(t))=f(u^{0}(t))$ if $t\in [0,\rT_{\infty}-\varepsilon ]$.
We also make the change of variable%
$$
\widetilde{t}=t-\tau
$$
and consider the delayed problem 
$$
\widetilde{\rP}(f,u^{0},\rB)=\left\{ 
\begin{array}{l}
\ry^{\prime }(t)=\lambda f_{\rM_{\epsilon}}(\ry)+\rB^{\prime }(t)\ry(t-\tau ),\quad 0<t<\tau 
\\ [.15cm]
\ry(\theta )=u^{0} (\theta ),\quad -\tau \leq \theta \leq 0
\end{array}
\right.
$$
(where, for simplicity, we have denoted again $\widetilde{t}$ by $t,$ so that, for
any $-\tau \leq \theta \leq 0$ we are identifying $u^{0} (\theta )$ with $
u^{0} (\theta+\rT_{\infty}-\epsilon )$, for some suitable function $\rB(t)$. Here $u^{0}(t)$ denotes again the solution without control.  Our goal is to
show that we can chose the control in the form
$$
\alpha(t):=\rB^{\prime }(t)\ry(t-\tau )
$$
such that the solution of $\widetilde{\rP}(f,u^{0},\rB)$ is defined on the whole
interval $[0,\tau )$ and that $\alpha\in \rW^{-1,q\prime }(0,\tau :\RR)$, for some $q>1$.

\par
Now, let us indicate the choice of function $\rB$ and the reformulation of $\widetilde{\rP}(f,u^{0},\rB)$ as a "neutral" equation. Given $q>1,~a>0,~\gamma \in (0,\frac{1}{q})$ and $
~ $a continuous function $m$  (to be taken in order to have $
\rB(0)=0 $, $\rB(t)>0$  and $\rB^{\prime }(t)>0$ on $0<t<\tau$) we define

$$
\rB(t)=\dfrac{a}{\left\vert t-t^{\ast }\right\vert ^{\gamma }}+m(t),\quad t\in [0,\tau ],
$$
with $t^{\ast }=\epsilon $ in this new time scale ($i.e.~ t=\rT_{\infty}$ in
the original time scale). We assume that $t^{\ast }\in (0,\tau ),~ i.e., ~
2\epsilon <\rT_{\infty}$. One possibility to avoid the difficulty related to the singularity of $\rB^{\prime }(t)$ is to reformulate $\widetilde{\rP}
(f,u^{0},\rB)$ (as in \cite{CDV2008}) as the "neutral" problem 
\begin{equation}
\left\{ 
\begin{array}{l}
\dfrac{d}{dt}\left[ \ry(t)-\rB(t)\ry(t-\tau )\right] 
=\lambda f_{M_{\epsilon}}(\ry)-\rB(t)\dfrac{d}{dt}\left[ \ry(t-\tau )\right], t\in [0,\tau ], \\ [.2cm]
\ry(\theta )=u^{0}(\theta ),\quad -\tau \leq \theta \leq 0.
\end{array}
\right. 
\label{Neutral ODE f}
\end{equation}

In addition, we will use the extension to the case of "neutral" equations of the version of the  {\em Alekseev's nonlinear variation of
constants formula } \cite{Alekseev} given in Theorem \ref{teo:theo2}). We recall that (for regular functions), this formula can be stated in the following terms:

\begin{prop} [Alekseev's formula, \cite{Alekseev}] Let $h:\RR\rightarrow \RR$ be $\cC^{2}$. Let 
$\ry^{0}(t)=\phi (t,t_{0},\xi )$ be the unique solution of the ODE 
$$
\left\{ 
\begin{array}{l}
\ry^{\prime }(t)=h(\ry(t)),\\ [.15cm]
\ry(t_{0})=\xi ,
\end{array}
\right.
$$
and let $\Phi (t,t_{0},\xi )=\partial _{\xi }\phi (t,t_{0},\xi ),$ where $
\partial _{\xi }$ denotes the partial differentiation. Then $\phi $ is $
\cC^{2},~\Phi $ is $\cC^{1}$, and for any $\rH:\RR\rightarrow \RR$ in $\rL_{
\text{loc}}^{1}$, the solution $z(t)$ of the so-called "perturbed" problem
$$
\left\{ 
\begin{array}{l}
z^{\prime }=h(z(t))+\rH(t),\\ [.15cm]
z(t_{0})=\xi ,
\end{array}
\right.
$$
has the integral representation 
$$
z(t)=\ry^{0}(t)+\int_{t_{0}}^{t}\Phi (t,s,z(s))\rH(s)ds.  \label{9}
$$
\fineq
\end{prop}
\begin{rem}\rm Notice that $\Phi (t,t_{0},\xi )$ satisfies $
\Phi (t,t,\xi )=1$.  Alekseev's
formula will be extended in Theorem \ref{teo:theo2} under a much greater generality and then applied in the framework of "neutral" equations.$\fin$
\end{rem}

Now, we can consider the
delayed term as an external "forcing"
$$
\rH(t)=\rB^{\prime }(t)\ry^{0} (t-\tau ),
$$
so that, by setting $t_{0}=0,~\xi =z(0)=u^{0}(0)=u_{0},~\ry^{0}(t)=\phi (t,0,\xi )$, the Alekseev's
formula  we can write (at least formally) 
$$
z(t)=\ry^{0}(t)+\int_{0}^{t}\Phi (t,s,z(s))\rB^{\prime }(s)u^{0}(s-\tau )ds.
$$
If we approximate $\rB(t)$ by regular functions (denoted again by $\rB(t)$), then the above formula can be equivalently written, after integrating by parts, as
$$
\begin{array}{lll}
z(t) & = & \disp \ry^{0}(t)+\bigg [ \Phi (t,s,z(s))\rB(s)u^{0}(s-\tau )\bigg ]
_{s=0}^{s=t}-\int_{0}^{t}\rB(s)\dfrac{d}{ds}\left[ \Phi (t,s,z(s))u^{0}(s-\tau )
\right] ds\\ [.275cm] 
& = & \disp \ry^{0}(t)+\Phi (t,t,z(t))\rB(t)u^{0}(t-\tau )-\int_{0}^{t}\rB(s)\dfrac{d}{ds}\left[ \Phi (t,s,z(s))u^{0}(s-\tau ) \right] ds.
\end{array}
$$
By the above remark, $\Phi (t,t,z(t))=1.$ On the other hand, as we saw
before, for $u^{0}\in \rW^{1,q}(-\tau ,0)$ its product by the $\cC^{1}$ function $\Phi (t,s,z(s))$ is also in $\rW^{1,q}(-\tau ,0)$. Therefore, its
derivative belongs to $\rL^{q}(-\tau ,0)$ and the indefinite integral, as in
all the previous cases, is an absolutely continuous function. Moreover the regularity of function $h$ is not needed in the final conclusion and thus we can argue by approximation (as we will make in the proof of Theorem 2). This means
that the integration by parts is legitimate and we may state the following
result, which is an extension of the Alekseev's
formula to "neutral" equations:
\begin{prop} The initial value problem
$$
\widetilde{\rP}(f,u^{0},\rB)=\left\{ 
\begin{array}{l}
\ry^{\prime }(t)=\lambda f_{\rM_{\epsilon}}(\ry) +\rB^{\prime }(t)\ry(t-\tau ),\quad 0<t<\tau 
\\ [.15cm]
\ry(\theta )=u^{0}(\theta ),\quad -\tau \leq \theta \leq 0%
\end{array}%
\right.
$$
\noindent with $f$ Lipschitz continuous and
initial function $u^{0}$ in $\rW^{1,q}(-\tau ,0)$ has a precise integral sense
in $[0,\tau ]$ by means of the neutral equivalent equation \eqref{Neutral ODE f}, and its unique
solution $z$ admits the integral representation 
\begin{equation}
z(t)=\ry^{0}(t)+\rB(t)u^{0}(t-\tau )-\int_{0}^{t}\rB(s)\dfrac{d}{ds}\big [ \Phi
(t,s,z(s))u^{0}(s-\tau )\big ] ds,  \label{representation z}
\end{equation}%
(where $\ry^{0}(t)=\phi (t,0,u^{0}(0)))$ associated to $h=\lambda f_{\rM_{\varepsilon}}$). Then, for every $u^{0}(\cdot) \in \rW^{1,r}(0,\tau
)$ (where $1/q+1/r=1)$ the neutral Cauchy problem has a unique solution
given by the identity \eqref{representation z}. 
\end{prop}

Therefore, $z\in \rL^{q}(0,\tau
),~z(t)-\rB(t)u^{0}(t-\tau )$ is an absolutely continuous function on $(0,\tau)$,
and then we may write 
$$
z(t)=\rB(t)u^{0}(t-\tau )+{\rm AC},
$$
where AC means an "absolutely continuous" function on the closed interval $[0,\tau
]$ (notice that without the truncation operation this is not necessarily true). As a consequence, the
singularity of the solution on $[0,\tau ]$ coincides with the singularity of $\rB$.
In particular, since $t^{\ast }=\epsilon $ (recall that $t^{\ast
}=\rT_{\infty}$ in the original scale of time), and by taking $0<\gamma <1$, we have that for some $m$  continuous function on $[0,\tau ]$ we have 
$$
\rB(t)=\dfrac{a}{\left\vert t-t^{\ast }\right\vert ^{\gamma }}+m(t).
\label{bb1}
$$
Since the initial function $u^{0}(\cdot)$ satisfies $u^{0}(t^{\ast }-\tau
)=\ry^{0}(\epsilon )\neq 0,$ then $t^{\ast }$ is also a singularity of $z$
(the controlled explosion) and 
$$
z(t)\simeq \dfrac{a}{\left\vert t-t^{\ast }\right\vert ^{\gamma }}
\ry^{0}(\epsilon ),\quad \text{as }t\rightarrow t^{\ast },
$$
is an asymptotic expansion of $z$ near $t^{\ast }=\rT_{\infty},$ which gives
the qualitative picture of the behavior of the solution near singularities
of $\rB$. Obviously, from the choice of $\gamma $ we get, finally that 
$z\in \rL^{q}(0,\tau)~$, which implies the desired property in this step: $\ry\in \rL^{1}\big (0,\rT_{\infty})$.

Moreover, the control $
\alpha(t)\doteq\rB^{\prime }(t)\ry(t-\tau )$ is in $\rW^{-1,q\prime }(0,\tau :\RR)$). Finally, notice that $u_{0}>0$ implies that $u^{0}(\theta )>0$ \ for\quad any $-\tau \leq \theta \leq 0$. Then, by
construction, we get that $u^{\alpha }(t)>0$ for any $t\in[0,\rT_{\infty}]$.

\par
\begin{rem}\rm \label{rem:bang-bang} The integrability condition $u^{\alpha }\in
\rL^{1}(0,\rT_{\infty})$ holds, in some special cases of the forcing
term $f(y)$ without any truncation and with $\alpha =0.$ This is the special case in which the function $\Phi
^{-1}\in \rL^{1}(0,\tau )$ for some $\tau >0.$ In the case of powers, $f(u)=u^{p}$, it
corresponds to the additional condition $p>2.$ Notice that this explains the non-uniqueness of the searched control. Some kind of optimality criterion on the set of searched controls could be introduced (see, e.g. the paper \cite{Amman-Quittner} and its references) but we will not enter here in this kind of considerations since the required techniques are of different nature. $\fin$
\end{rem}
\par
\noindent {\em Step 2}: $t\in [\rT_{\infty},2\rT_{\infty})$. We will take a control such that $\alpha (t)\leq 0$ on $(\rT_{\infty},2\rT_{\infty}).$ Thus, we consider the problem    
\begin{equation}
\left\{ 
\begin{array}{lc}
\dfrac{du}{dt}(t)+\lambda f_{\rM_{\varepsilon}}(u)=\alpha (t) & \text{in }%
(\rT_{\infty},2\rT_{\infty}), \\ [.175cm]
u(\rT_{\infty})=+\infty . & 
\end{array}%
\right.   
\label{Eq: second iteration}
\end{equation}
As mentioned before, the nonlinear term is of absorption type and problems of this nature were already considered in  \cite{BaDiDi} (see problem (D)) but without any truncation argument and with  $\alpha =0$. The case with a truncation and $\alpha \leq 0$ can be solved by reflection with respect to the time $\rT_{\infty}$. Indeed, we define
$$
\rY_{\widehat{\alpha }}(t)=u^{\alpha }(t-\rT_{\infty})\text{ for }t\in [
\rT_{\infty},2\rT_{\infty}]\quad \text{ and }\quad \widehat{\alpha }(t)=-\alpha (t-\rT_{\infty}).
$$
It is a routine matter to check that $Y_{\widehat{\alpha }}(t)$ satisfies
problem \eqref{Eq: second iteration} for the control $\widehat{\alpha }(t)$. 
Notice that, $\widehat{\alpha }(t)=0$ on the interval $(\rT_{\infty}+%
\varepsilon ,2\rT_{\infty}]$, $\widehat{\alpha }(t)<0$ on ($%
\rT_{\infty},\rT_{\infty}+\varepsilon )$ (in fact  $\widehat{\alpha }(t)\searrow
-\infty $, if $t\searrow \rT_{\infty}),~\rY_{\widehat{\alpha }}(t)>0$ on $%
(\rT_{\infty},2\rT_{\infty}],~\rY_{\widehat{\alpha }}(2\rT_{\infty})=u_{0}$ and $\rY_{
\widehat{\alpha }}\in\rL^{1}(\rT_{\infty},2\rT_{\infty}).$
\par
\noindent {\em Step 3}. Finally, we use a $2\rT_{\infty}$-periodicity argument to extend the controlled solution to the interval $t\in (2\rT_{\infty},+\infty)$. It is clear that for times in which the control changes sign we get some singularity on the time derivative of the solution but, as typical in control theory, the differential equation holds for almost any $t\in (0,+\infty)$, the controlled solution is in $\rL_{loc}^{1}(0,+\infty)$, and the proof of Theorem 1 is complete.
\begin{figure}[htp]
	\begin{center}
		\includegraphics[width=15cm]{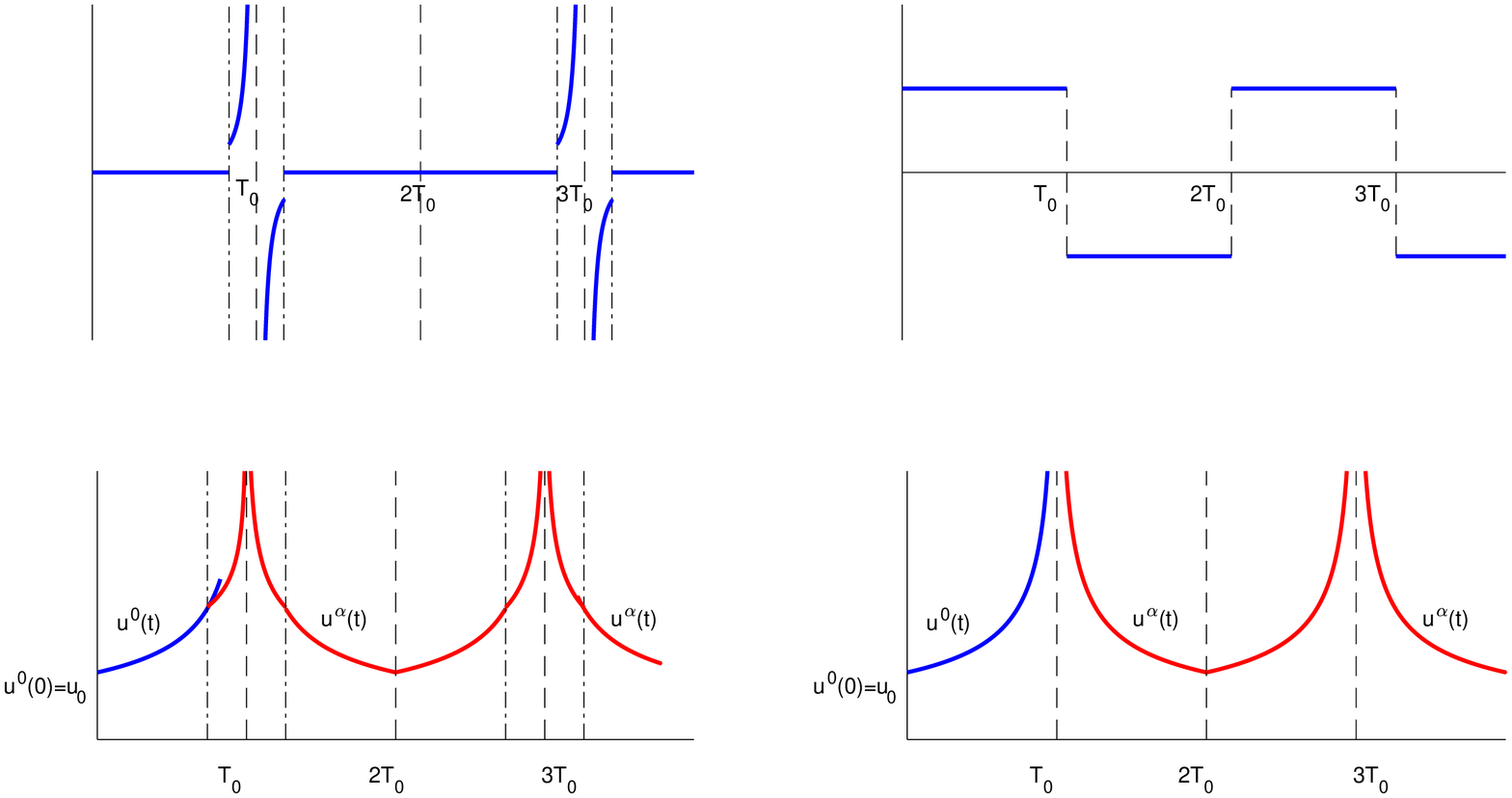}\\ [-.5cm]
		\begin{tabular}{l}\small Figure 1. Illustrative example of the control $\alpha(t)$ and the effective bang-bang control if $f(s)=s^{p},p>2$.\\ \small
		The blowing up solution without control $u^{0}(t)$ and the controlled solution $u^{\alpha}(t)$ defined in the whole $[0,+\infty[$.
	    \end{tabular}
	\end{center}
	\label{figure:controls}
\end{figure}
\begin{exam}\rm We consider the special case of $%
f(u)=u^{2}$ and $u_{0}=1$. Then we can identify easily the elements
appearing in the proof of Theorem 1. Indeed, in this case, 
$$
\phi (t,t_{0},\xi )=\frac{1}{\frac{1}{\xi }-(t-t_{0})},
$$
and so $\rT_{\infty}=1$. Thus we can take, $e.g.,~\varepsilon =1/8$ (and, of
course, $2\epsilon <\rT_{\infty}$), $\tau =\rT_{\infty}-\varepsilon =7/8$. Taking $%
\gamma =1/5$ and $a=1$ we get that  $\rB(t)=\dfrac{1}{\left\vert 1-t\right\vert
^{5}}$ and thus the searched control $\alpha (t)$ is given by%
$$
\alpha (t)=\left\{ 
\begin{array}{cc}
0 & \text{if }t\in (0,7/8) \\ [.15cm]
\rB^{\prime }(t)u(t-7/8) & \text{if }t\in (7/8,1),%
\end{array}%
\right. 
$$
with $u(t)$ the solution of the delayed problem%
$$
\left\{ 
\begin{array}{l}
u^{\prime }(\widetilde{t})=f_{\rM_{\varepsilon }}(u(\widetilde{t}))+\dfrac{5}{
\left\vert 1-\widetilde{t}\right\vert ^{6/5}}(u(\widetilde{t}-7/8)),\quad 
\widetilde{t}\in (7/8,1)\\ [.45cm]
u(\theta )=u^{0}(\theta ),\quad -7/8\leq \theta \leq 0,
\end{array}%
\right.   
$$
where $u^{0}(\theta )=\dfrac{1}{1-\theta }$ if $\theta \in  [-7/8,0]$ and 
$$
f_{\rM_{\varepsilon }}(u)=\left\{ 
\begin{array}{cc}
u^{2} & \text{if }u\in (0,7/8) \\ [.15cm]
49/64 & \text{if }u\in (7/8,+\infty ).
\end{array}%
\right. 
$$
\fineqnum
\end{exam}

\subsection{Proof of Theorem 2}
\noindent {\sc Proof of Theorem \ref{teo:theo2}}. Let $h_{n}\in \cC^{1}(\RR^{d}:\RR^{d})$ be a sequence approximating $h$ in $\rW^{1,s}(\RR^{d}:\RR^{d})$, for any $s\in [1,+\infty ),$ and such that 
\begin{equation}
\left\Vert \partial _{x}h_{n}\right\Vert
_{\rL^{\infty }(\RR^{d}:\cM_{d\times d})}\leq \left\Vert
\partial _{x}h\right\Vert _{\rL^{\infty }(\RR^{d}:\cM_{d\times d})}\doteq\rM\text{ for any }n\in \NN\hbox{}  \label{Lipschitz constant}
\end{equation}
Let $\ry_{n}^{0}=\phi
_{n}(t,t_{0},\xi )$ be the unique solution of the unperturbed ODE 
$$
\rP^{\ast }(h_{n},0,\xi )=\left\{ 
\begin{array}{l}
\ry^{\prime }(t)=h_{n}(\ry(t))\text{ \ in }\RR^{d}, \\ [.15cm]
\ry(t_{0})=\xi ,%
\end{array}
\right.
$$
and let $\Phi _{n}(t,t_{0},\xi )=\partial _{\xi }\phi _{n}(t,t_{0},\xi ),$.
Let us consider the sequence of "perturbed" problems 
$$
\rP^{\ast }(h_{n},\beta ,\xi )=\left\{ 
\begin{array}{l}
\dfrac{d\ry_{n}}{dt}(t)+\beta (t,\ry_{n}(t))\ni  h_{n}(\ry_{n}(t)),\text{ in }%
\RR^{d}, \\ 
\ry(t_{0})=\xi .%
\end{array}%
\right.
$$
Then, by the classical version of the Alekseev formula (also valid for $%
d\geq1$) we know that 
\begin{equation}
\ry_{n}(t)=\ry_{n}^{0}(t)-\int_{t_{0}}^{t}\Phi _{n}(t,s,\ry_{n}(s))\beta
(s,\ry_{n}(s))ds,\text{ for any }t\in [0,\rT],
\label{n-representation formula}
\end{equation}%
(as before, in the above formula we assumed, for simplicity, that $\beta (t,\cdot)$ is single-valued but a suitable similar expression
can be obtained if $\beta (t,$\textperiodcentered $)$ is multivalued). Since $h_{n}\rightarrow h$ \ and $h$ is locally Lipschitz, we know that $
\ry_{n}^{0}\rightarrow \ry^{0}$ and $\ry_{n}\rightarrow \ry$ strongly in ${\rm AC}([0,\rT]:\RR^{d})$ for any fixed $\rT>0$ (this is an easy application of Theorem 4.2 of Brezis \cite{BreOMM} for the autonomous case and from \cite{Ya} in the non-autonomous case). Moreover since any maximal monotone operator is strongly-weakly closed we know that, at
least, $\beta \big (\ry_{n}(\cdot),\cdot\big )\rightharpoonup \beta \big (\ry(\cdot),\cdot\big )$ in $\rL^{2}(0,\rT:\RR^{d})$. Then, from the classical Peano
Theorem we know that there exists a $\Phi (t,s,\ry)$ such that 
$$
\Phi _{n}\big (t,\cdot,\ry_{n}(\cdot)\big )\rightarrow \Phi \big (t,\cdot,\ry(\cdot)\big )\text{, for a.e. }t\in (0,\rT),
$$
strongly in $\rL^{2}(0,\rT:\cM_{d\times d})$. Indeed, $\Phi
_{n}(t,t_{0},\xi )$ is the solution of the problem 
$$
\left\{ 
\begin{array}{l}
\Phi ^{\prime }(t)=\rH_{n}(t,t_{0},\xi )\Phi (t)\text{ \ in }\cM
_{d\times d}, \\ [.15cm]
\Phi (t_{0})=\rI,
\end{array}
\right.
$$
where 
$$
\rH_{n}(t,t_{0},\xi )\doteq\partial _{x}h_{n}(\phi _{n}(t,t_{0},\xi )).
$$
But, we know that, if $\rM$ is given by \eqref{Lipschitz constant} then 
$$
\left\Vert \rH_{n}(t,t_{0},\xi )\right\Vert _{\rL^{\infty }(t_{0},\rT:\cM
_{d\times d})}\text{ }\leq \rM\text{ \ for any }t_{0}\in (0,\rT)\text{ and for
any }\xi \in \RR^{d}.
$$
Thus, by Gronwall inequality, there exists a positive constant $\widetilde{\rM}%
=\widetilde{\rM}(t_{0},\xi )$ such that
$$
\left\Vert \Phi _{n}(\cdot,t_{0},\xi )\right\Vert
_{\rW^{1,\infty }(0,\rT)}\leq \widetilde{\rM}.
$$
This implies that there exists a Lipschitz function $\Phi (t,s,\xi )$ such
that $\Phi _{n}\big (t,\cdot,y_{n}(\cdot%
)\big )\rightharpoonup \Phi \big (t,\cdot,y(\cdot)\cdot )$
in $\rW^{1,q}(0,\rT:\cM_{d\times d})$ for any $q\in (1,\infty ).$ This
leads to the strong convergence in $\rL^{2}(0,\rT:\cM_{d\times d})$.
Then we can pass to the limit in formula \eqref{n-representation formula}
and we get that 
$$
\ry(t)=\ry^{0}(t)-\int_{t_{0}}^{t}\Phi (t,s,\ry(s))\beta (s,\ry(s))ds,\text{ for any
}t\in [0,\rT].
$$
\fineqnum

\section{The complete recuperation after the blow up time for problem $\rP(\alpha)$}
\label{section:completerecuperation}
Previously to the consideration of the controlled problem, it is useful to establish some basic properties for the uncontrolled problem $\rP(0)$ 
$$
\left\{
\begin{array}{ll}
-\Delta u^{0}+g(u^{0})=0 & \quad\hbox{in }\bB_{\rR}\times(0,\infty),\\[.15cm]%
\dfrac{\partial u^{0}}{\partial t}+\dfrac{\partial u^{0}}{\partial\bn }=\lambda f(u^{0}) &
\quad\hbox{on }\partial\bB_{\rR}\times (0,\infty),\\[.25cm]%
u^{0}(\rR,0)=u_{0}>0, & \quad x\in\partial \bB_{\rR},
\end{array}
\right.  
$$
where $g$ and $f$ are continuous non-negative increasing real functions as indicated in the Introduction and $\lambda$ is a positive constant. 
The framework of our study concerns the case in which the forcing term dominates over the absorption one, in the sense of the condition \eqref{eq:dominationbalance}
$$
\liminf_{\tau\rightarrow\infty}\dfrac{f(\tau)}{\sqrt{2\rG\big (\tau\big )}}\in (0,+\infty ].
$$
The main goal of this section is not only to prove that it is possible to
build a control $\alpha (t)$ such that the solution $u^{\alpha }$ of
problem $\rP(\alpha)$ be well defined for any $t>0$,  but to
prove previously that in the isolated times in which the solution blows up it takes
place only on the boundary $\partial \rB_{\rR}$.  So, we will prove that at the blow up time the solution without control will coincide with the unique large solution of the problem 
\begin{equation}
	\left\{
	\begin{array}{ll}
		-\Delta\rU_{\infty}^{\bB_{\rR}}+g\big (\rU_{\infty}^{\bB_{\rR}}\big )=0 & \quad
		\hbox{in ${\bB_{\rR}}$},\\[.15cm]
		\rU_{\infty}^{\bB_{\rR}}=\infty & \quad\hbox{on $\partial {\bB_{\rR}}$.}
	\end{array}
	\right.  \label{eq:largesolution ball}
\end{equation}
Beside condition \eqref{eq:KellerOssermanG} we will require other technical assumptions:
\begin{equation}
	\limsup_{s\rightarrow\infty} \frac{\Psi(\eta s)}{\Psi(s)} <1
	\quad\hbox{for any } \eta>1 \label{eq:nopotencias}
\end{equation}
(see \eqref{eq:inverseexplosiveprofile} below) and 
\begin{equation}
	\dfrac{g(s)}{s}\quad\hbox{ is increasing for large $s$.}
	\label{eq:monotonicityfracg}
\end{equation}
Among other examples, the conditions \eqref{eq:KellerOssermanG}, \eqref{eq:nopotencias} and \eqref{eq:monotonicityfracg}
hold for the power-like case $g_{m}(s)=s^{m}$, when $m>1$.

It is clear that there are two different subcases in which the the forcing term dominates over the absorption one: 
\par
\noindent a) Forcing term strongly dominating over absorption term.
It corresponds to the case in which the following condition holds 
\begin{equation}
	\lim_{\tau\rightarrow\infty}\dfrac{f(\tau)}{\sqrt{2\rG\big (\tau\big )}}=\infty,
	\label{eq:fsuperwinG}
\end{equation}
with $\displaystyle \rG(s)=\int^{s}_{0}g(s)ds$. Here, for any $\lambda>0 $ the domination at infinity of the forcing term over the expression $\sqrt{2\rG}$ associated to the absorption  term is satisfied. 
\par
\noindent b) Forcing term weakly dominating over absorption terms.
It concerns the case in which we have
\begin{equation}
	\liminf_{\tau\rightarrow\infty}\dfrac{f(\tau)}{\sqrt{2\rG(\tau)}}=\rL>0. 
	\label{eq:fwinG-}
\end{equation}
In some sense, the functions $\sqrt{2\rG(\tau)}$ and $f(\tau)$ are of the same order at infinity. We will see that the domination at infinity of the effective forcing term  $\lambda f$ over $\sqrt{2\rG}$ requires the assumption 
$$
\lambda> \dfrac{1}{\rL},
$$
which is obvious when the forcing term strongly dominates to the absorption (condition \eqref{eq:fsuperwinG}).
\begin{rem}\rm For the power-like case $f_{p}(r)=s^{p},~p>0,$ and $g_{m}(s)=s^{m},~m>0$. It implies $\rG_{m}(s)=\dfrac{1}{m+1}s^{m+1}$. They  imply
$$
\dfrac{f_{p}(\tau)}{\sqrt{2\rG_{m}(\tau)}}=\sqrt{\dfrac{m+1}{2}}\tau^{\frac{2p-(m+1)}{2}}.
$$
Then, the condition \eqref{eq:fsuperwinG} hods if $2p>m+1$. On the other hand, \eqref{eq:fwinG-} requires the equality 
$2p=m+1$ and then $\rL=\sqrt{\dfrac{m+1}{2}}$ (we will come back on it in Theorem \ref{theo:invariance} below).$\fin$ 
\label{rem:powerdomination1}
\end{rem}
\par
In fact, we are interested in a global domination
\begin{equation}
	\lambda f(\tau)>\sqrt{2\rG(\tau)},\quad \hbox{for any }\tau>0,
	\label{eq:globaldomination}
\end{equation}
for $\lambda >0$. In some cases,  the local domination \eqref{eq:dominationbalance} can imply  the global domination for suitable large values of the parameter $\lambda$ (see (\ref{Dato inicial grande}).
\begin{prop}
Assume \eqref{eq:dominationbalance}. If   
\begin{equation}
	\liminf_{\tau\searrow 0}\dfrac{f(\tau)}{\sqrt{2\rG(\tau)}}>0,
	\label{eq:initialdomination}
\end{equation}
there exists $\lambda_{0}>0$, depending only on $f$ and $g$, for which one has the restricted global domination
\begin{equation}
	\lambda f(\tau)>\sqrt{2\rG(\tau)}\quad \hbox{for any }\tau>0,
	\label{eq:globaldominationobtenida}
\end{equation}
provided $\lambda >\lambda _{0}$. Therefore a global domination is verified under \eqref{eq:globaldomination} or under the couple of condition 
\eqref{eq:dominationbalance} and \eqref{eq:initialdomination}.
\par
\noindent 
Whenever \eqref{eq:initialdomination} fails, thus
$$
\liminf_{\tau\searrow 0}\dfrac{f(\tau)}{\sqrt{2\rG(\tau)}}=0,
$$
we only may extend the local domination \eqref{eq:dominationbalance}  to
$$
	\lambda f(\tau)>\sqrt{2\rG(\tau)},\quad \tau>\tau_{-},\quad \hbox{for each $\tau_{-}>0$}
$$
for  $\lambda >\lambda _{0}$. Here $\lambda_{0}>0$, depending on $f, g$ and $\tau_{-}$.
\label{prop:wideregiondomination}
\end{prop}
\proof  The assumption \eqref{eq:fwinG-} implies that there exists $\tau _{0}>0$, large enough, such that
\begin{equation}
	\dfrac{f(\tau)}{\sqrt{2\rG(\tau)}}>\rL>0 \quad \hbox{for any }\tau >\tau_{0}.
	\label{eq:minorationdomination}
\end{equation}
Since $f$ and $\sqrt{2\rG}$ are positive continuous  functions, we deduce
$$
\dfrac{f(\tau)}{\sqrt{2\rG(\tau)}}>\rL_{*}\doteq \min_{0<\tau	\le\tau_{1}}\dfrac{f(\tau)}{\sqrt{2\rG(\tau)}}\ge 0,\quad \hbox{for }0<\tau\le \tau_{0}. 
$$
Therefore \eqref{eq:initialdomination} implies
$$
\dfrac{f(\tau)}{\sqrt{2\rG(\tau)}}>\rL_{0}\doteq \min\{\rL,\rL_{*}\}>0, \quad \hbox{for any }\tau >0,
$$
and then \eqref{eq:globaldominationobtenida} holds for $\lambda >\lambda_{0}\doteq \dfrac{1}{\rL_{0}}$. The reasoning also applies to the case in which \eqref{eq:fsuperwinG} holds and  we have an inequality similar to  \eqref{eq:minorationdomination} for any $\rL>0.\fin$
\begin{rem}\rm Since for  $f_{p}(s)=s^{p}$ and $g_{m}(s)=s^{m}$ one has
$$
\dfrac{f_{p}(\tau)}{\sqrt{2\rG_{m}(\tau)}}=\sqrt{\dfrac{m+1}{2}}\tau^{\frac{2p-(m+1)}{2}},
$$
the condition \eqref{eq:initialdomination} holds if $2p= m+1$ (see Theorem \ref{theo:invariance}).$\fin$
\label{rem:powerdomination}
\end{rem}
\medskip

\par
The main result of this Section is the following: 
\begin{theo}\label{Thm PIII} 
Assume \eqref{eq:KellerOssermanG}, \eqref{eq:dominationbalance}, \eqref{Dato inicial grande} and \eqref{eq:globaldomination}. Then:

\noindent i) For any $u^{0}(\rR,0)=u_{0}>0$  there exists a finite time $\rT_{\infty} (u^{0})$ such that the solution without control coincides at this time with the large solution of the problem \eqref{eq:largesolution ball}. In fact the continuation, $u^{0}(\cdot,t)=\rU_{\infty }^{\rB_{\rR}}(\cdot)$ for $t>\rT_{\infty}(u^{0})$ is a global solution of problem $\rP(0)$ blowing up on the boundary for any $t>\rT_{\infty}(u^{0})$.

\noindent ii) For any $u_{0}>0$ large enough, the blowing up trajectory $u^{0}(\cdot,t)$  of the associated problem $\rP(0)$ has a controlled
explosion {\rm (}in the sense of Definition 1{\rm )} by means of the control problem $\rP(\alpha)$, with the structure conditions (\eqref{hyp purebangnbang}) and (\eqref{hyp quasibangnbangGNew}), for a suitable $\alpha \in \rW_{loc}^{-1,q'}(0,\rT_{\infty} :\RR)$, for some $ q>1.$ Moreover, if $u^{\alpha }(\cdot,t)$, in $\rL_{loc}^{1}(0,+\infty :\rL^{\infty}(\Omega))$, is the solution corresponding to the built control $
\alpha (t)$, then $u^{\alpha }(x,t)<+\infty$ for any $t>0$ and any $x$ in $\bB_{\rR}$.
\end{theo}

As above we use the notation $ \rT_{\infty}(u^{0})$ since this time depends not only on the initial datum $u_{0}$ but the other parameters. Again in the following we will simplify the notation by writing  $\rT_{\infty}(u^{0})=\rT_{\infty}$.

The crucial step in our study is to consider previously the case without any control $\rP(0)$. We will need
some previous results (collected in Section 3.1) and to get the existence
and uniqueness of solutions of problem $\rP(0)$ and to prove part i) of Theorem \ref{Thm PIII}. Finally, in Section 3.3 we will give the proof of part ii) of Theorem \ref{Thm PIII}.

\subsection{The boundary blow up for the uncontrolled problem}
\label{section:profile} 

Some of the results of this Section are applicable to the case of an arbitrary (non necessarily symmetric)  open set $\Omega\subset \RR^{\rN},~\rN>1$.  We collect here some useful technical properties derived from the condition \eqref{eq:KellerOssermanG}. 
\begin{lemma}
[Lemma 6.1 of \cite{AlDiRe}]Let us assume \eqref{eq:KellerOssermanG}. Then
\begin{equation}
\lim_{s\rightarrow\infty} \dfrac{s}{\sqrt{\rG(s)}}
=\lim_{s\rightarrow\infty}\dfrac{s}{g(s)}=\lim
_{s\rightarrow\infty} \frac{\sqrt{\rG(s)}}{g(s)}=0,
\label{eq:potenciamayorqueuno}
\end{equation}
where $\displaystyle \rG(s)=\int^{s}_{0}g(s)ds$. Hence
$$
\dfrac{s}{g(s)}=o\left(  \dfrac{s}{\sqrt{\rG(s)}
}\right)  \quad\hbox{and}\quad\dfrac{s}{g(s)}=o\left(  \frac
{\sqrt{\rG(s)}}{g(s)}\right).
$$
\fineq
\label{Lemma:KellerOssermanAGfraccion}
\end{lemma}
\par
An useful tool in this subsection is the decreasing function associated to the improper finite integral \eqref{eq:KellerOssermanG}, 
\begin{equation}
\Psi(\delta)=\int^{\infty}_{\delta}\dfrac{ds}{\sqrt{2\rG(s)}}
\quad\hbox{for any}\delta>0 
\label{eq:inverseexplosiveprofile}
\end{equation}
with $\Psi(\infty)=0$. Straightforward computations allows to see that for $\zeta $ small we have 
\begin{equation}
\dfrac{d }{d \zeta}\Psi^{-1}(\zeta)=-\sqrt{2\rG\big (\Psi^{-1}(\zeta)\big )}\quad \hbox{and}
\quad \dfrac{d ^{2}}{d \zeta^{2}}\Psi^{-1}(\zeta)=g\big (\Psi^{-1}(\zeta)\big ).
\label{eq:tecnnivalityPsi}
\end{equation}
Properties \eqref{eq:tecnnivalityPsi} can be used to characterize the unique explosive
profile on the boundary of the large solution of the problem
\begin{equation}
\left\{
\begin{array}{ll}
-\Delta\rU_{\infty}^{\Omega}+g\big (\rU_{\infty}^{\Omega}\big )=0 & \quad
\hbox{in $\Omega$},\\[.15cm]
\rU_{\infty}^{\Omega}=\infty & \quad\hbox{on $\partial \Omega$,}
\end{array}
\right.  \label{eq:largesolution}
\end{equation}
provided \eqref{eq:monotonicityfracg}.
More precisely
\begin{theo}{{\rm \cite{DiLe,AlDiRe}}}
Let $\Omega\subset\RR^{\rN},~\rN>1$ be a bounded open set
where $\partial\Omega$ satisfies an inner and outer sphere condition. Assume
\eqref{eq:KellerOssermanG}, \eqref{eq:nopotencias} and  \eqref{eq:monotonicityfracg}.
Then there exists a unique classical solution $\rU_{\infty}^{\Omega}$ of  \eqref{eq:largesolution} whose explosive boundary profile satisfies 
\begin{equation}
\lim_{\hbox{\rm dist($x,\partial \Omega$)}\rightarrow0}\dfrac{\rU_{\infty
	}^{\Omega}(x)}{\Psi^{-1}\big (\hbox{\rm dist($x,\partial \Omega$))}}=1.
\label{eq:boundaryprofile}
\end{equation}
\fineqnum
\label{theo:blowup}
\end{theo}
In fact, when $\Omega=\bB_{\rR}$ the rate behaviour \eqref{eq:boundaryprofile} becomes 
$$
\lim_{|x|\nearrow \rR}\dfrac{\rU_{\infty
}^{\rB_{\rR}}(x)}{\Psi^{-1}\big (\rR-|x|)}=1.
$$
\begin{rem}\rm
For the power like case $g(s)=s^{m}$, condition \eqref{eq:KellerOssermanG} becomes $m>1$ and
$$
\Psi_{m}(\delta)=\dfrac{\sqrt{2(m+1)}}{m-1}\dfrac{1}{\delta^{\frac{m-1}{2}}}
,\quad\delta\ge0.
$$
Moreover, the technical conditions \eqref{eq:nopotencias} and  \eqref{eq:monotonicityfracg} also hold. Then
$$
\rU_{\infty}^{\Omega}(x)=\left(  \frac{2(m+1)}{(m-1)^{2}}\right)  ^{\frac{1}{m-1}
}\big (\hbox{dist($x,\partial \Omega$)}\big )^{-\frac{2}{m-1}}
+o\big (\hbox{dist($x,\partial \Omega$)}\big )
$$
(see \cite{AlDiRe}).$\fin$ 
\label{rem:potenciasg}
\end{rem}

\begin{rem}\rm
Among other illustrative choices satisfying \eqref{eq:KellerOssermanG},
\eqref{eq:nopotencias} and \eqref{eq:monotonicityfracg} studied in \cite{AlDiRe} we pick up the function $g(s)=e^{s}$ for which 
$$
\rU_{\infty}^{\Omega}(x)=\log \left (\dfrac{2}{\big (\hbox{dist($x,\partial \Omega$)}\big )^{2}}\right )+o\big (\hbox{dist($x,\partial \Omega$)}\big ).
$$
Another example is $g(s)=se^{2s},$ for which
$$
\rU_{\infty}^{\Omega}(x)=\sqrt{2}\hbox{ erfc}^{-1}\left(\dfrac{\hbox{dist($x,\partial \Omega$)}}{\sqrt{\pi}}
\right)+o\big (\hbox{dist($x,\partial \Omega$)}\big ),
$$
where $\disp \hbox{erfc}(\delta)=1-\hbox{erf}(\delta)=\dfrac{2}
{\sqrt{\pi}}\int^{\infty}_{\delta}e^{-s^{2}}ds.\hfill\quad_{\Box}$
\label{rem:exponential} 
\end{rem}
\begin{rem}\rm Reasoning as in Remark \ref{rem:exemplesf}, from the inequality
$$
\widehat{g}(s)\ge g(s)\quad \hbox{for large $s$},
$$
it follows that if $g$ verifies  \eqref{eq:KellerOssermanG} the same happens with $\widehat{g}$. In particular, any function $\widehat{g}(s)\ge s q(s)$ verifying 
$$
\liminf_{s \rightarrow \infty} \dfrac{q(s)}{s^{\gamma}}\in (0,+\infty]\quad 
\hbox{for some $\gamma >0$}
$$
satisfies \eqref{eq:KellerOssermanG}. For instance, we may choose  $q(s)\ge (\log s)^{\gamma},~\gamma\ge 1,$ or $q(s)\ge \log (\log (\cdots \log (s)\cdots ))$.
It is also clear that an assumption like
$$
\dfrac{g(s)}{s^{\alpha}}\quad\hbox{ increasing for large $s$},
$$
for some $\alpha >1$, implies \eqref{eq:KellerOssermanG}.$\fin$
\label{rem:exemplesg}
\end{rem}
\begin{rem}\rm Sometimes it is more useful to write \eqref{eq:nopotencias} as 
$$
\liminf_{s\rightarrow\infty} \frac{\Psi(\eta s)}{\Psi(s)}>1
\quad\hbox{for $ 0<\eta<1$}.
$$
We note that for some example as $g(s)=s\big (\log s)^{m},~m>2$, which verifies \eqref{eq:KellerOssermanG}, the condition \eqref{eq:nopotencias} fails. In \cite{AlDiRe} a sharp argument enables us to extend Theorem \ref{theo:blowup} to the so called borderline case given by
\begin{equation}
\limsup_{s\rightarrow\infty} \frac{\Psi(\eta_{0} s)}{\Psi(s)}=1
\quad\hbox{for some $\eta_{0}>1$}. \label{eq:nopotenciasborder}
\end{equation}
This happens for the above choice or when $g(s)=c_{1}s(\log s)^{m}+c_{2}(\log s)^{m-1}
,~s>1,~m>2$, for $c_{1}>0$ and $c_{2}\in\RR$  or $g(s)=s(\log s)^{2}\big (\log (\log (s))\big )^{m},~s>0,~m>2$. We send to  \cite{AlDiRe} for some comments and other examples. We will not consider the borderline case in this paper.$\fin$
\label{rem:nopotencias}
\end{rem}
We note that the spatial explosive profile given by \eqref{eq:boundaryprofile} does
not depend on the geometrical properties of $\rB_{\rR}$ as
curvature or dimension. These influences can appear in lower term of the
explosive expansion near $\partial\rB_{\rR}$ (see again~\cite{AlDiRe}). We point out that 
some authors have approached the boundary behaviour of the large solutions, $e.g.$ Bandle, Ess\'en, Lazer, Marcus, Mc Kenna, Matero and many others (see \cite{Bandle, DiLe,AlDiRe}
and the references therein).

\subsection{Blowing up time-profile for the uncontrolled problem: the radially symmetric case}

Here we  will get some growing time-estimates near the blow-up time for problem  $\rP(0)$  on the radial domain $\bB_{\rR}$. The radially symmetric solution $u^{0}(x,t)\doteq u^{0}(|x|,t)$, corresponding to a constant initial datum satisfy
\begin{equation}
\left \{
\begin{array}{ll}
-\dfrac{1}{r^{\rN-1}}\dfrac{\partial u^{0}}{\partial r}\left( r^{\rN-1}\dfrac{\partial u^{0}}{\partial r}(r,t)\right)  +g\big (u^{0}(r,t)\big )=0,& \quad r<\rR,~ t>0\\[.35cm]
\dfrac{\partial u^{0}}{\partial r}(0,t)=0, &\quad t\ge0,\\ [.35cm]
\dfrac{\partial u^{0}}{\partial r}(\rR,t)+\dfrac{\partial u^{0}}{\partial r }(\rR,t)=\lambda f\big (u^{0}(\rR,t)\big ), &\quad t\ge 0,\\[.3cm]
u^{0}(\rR,0)=u_{0}>0. &
\end{array}
\right.  \label{eq:uncontrolledradialproblem}
\end{equation}
As the center of the ball does not play any important role we may assume that it is the origin of the space.  
\begin{theo}[ Blow up time on the boundary]  Assume  $\lambda >0$ satisfying the global domination of the Proposition \ref{prop:wideregiondomination}. Assume also  \eqref{eq:KellerOssermanG}. 
Then $\rP(0)$ has a unique radially symmetric solution, $u^{0}(|x|,t)$, on $\overline{\rB}_{\rR}\times [0,\rT_{\infty}\big )[$, for some  $\rT_{\infty}\le \Psi (u_{0})$, such that
$$
\left \{
\begin{array}{l}
0\le u^{0}(|x|,t)<\rU_{\infty}^{\bB_{\rR}}(x),\quad  (x,t)\in \bB_{\rR}\times [0,\rT_{\infty}[,\\ [.15cm]
\disp \lim_{t\nearrow \rT_{\infty}}u^{0}(|x|,t)=\rU_{\infty}^{\bB_{\rR}}(x),\quad x\in \overline{\bB}_{\rR},
\end{array}
\right .
$$ 
where $\rU_{\infty}^{\bB_{\rR}}$ is the relative stationary large solution on the ball $\bB_{\rR}$ {\rm (}see \eqref{eq:largesolution}{\rm )}. 
Moreover, under \eqref{eq:nopotencias}  the solution of \eqref{eq:uncontrolledradialproblem} has the explosive boundary behaviour
$$
\liminf_{t\nearrow \rT_{\infty}(u_{0})}\dfrac{u^{0}(\rR,t)}{\Psi^{-1}\big (\rT_{\infty}-t\big )\big)}\ge 1.
$$
\label{theo:boundaryblowup} 
\end{theo}
\begin{rem}\rm As it was deduced from the below proofs,  in some case as for the power choices  $f_{p}(s)=s^{p}$ and $g_{m}(s)=s^{m}$ with $2p>m+1$ the constant $\lambda_{0}$, introduced in Proposition \ref{prop:wideregiondomination} can also depend on the data $\rR$ and $u_{0}.\fin$   
\label{rem:lambda0RT}
\end{rem}
The proof of Theorem \ref{theo:boundaryblowup} is based on a result taking the advantage that for radially symmetric functions the comparison principle holds even if there are some singularities  at the origin.

\begin{prop}
Assume  $\lambda >0$ satisfying the global domination of the Proposition \ref{prop:wideregiondomination} as well as \eqref{eq:KellerOssermanG}. 
Define the function 
\begin{equation}
\underline{\rU}(r,t)=\Psi^{-1}\big (\nu(\rT-t +\rR-r\big )\big), \quad 0\le r<\rR,~0\le t<\rT,
\label{eq:explosivesubsolution}
\end{equation}
for $\nu >1$ and $\rT>0$. Then we have
$$
\left\{
\begin{array}{ll}
-\dfrac{1}{r^{\rN-1}}\dfrac{\partial \underline{\rU}(r,t)}{\partial r}\left( r^{\rN-1}\dfrac{\partial \underline{\rU}(r,t)}{\partial r}\right) +g\big (\underline{\rU}(r,t)\big )=\underline{\rE}(r,t),& \quad r<\rR,~ 0<t<\rT,\\[.35cm]
\underline{\rU}(0,t)=\Psi^{-1}\big (\nu(\rT-t +\rR \big )\big), &\quad 0\le t<\rT,\\ [.35cm]
\dfrac{\partial \underline{\rU}(0,t)}{\partial r}=\nu\sqrt{2\rG\big ( \underline{\rU}(0,t)\big )}, &\quad 0\le t<\rT,\\ [.35cm]
\dfrac{\partial \underline{\rU}(\rR,t)}{\partial r}+\dfrac{\partial \underline{\rU}(\rR,t)}{\partial r }\le \lambda f\big (\underline{\rU}(\rR,t)\big ), &\quad 0\le t<\rT,\\[.35cm]
\underline{\rU}(\rR,0)=\Psi^{-1}\big (\nu\rT\big), &\\
\end{array}
\right.
$$
where
\begin{equation}
\underline{\rE}\in \cC\big ((0,\rT):\rL^{2}(\varepsilon,\rR)\big ),\quad \hbox{for any $\varepsilon>0$}
\label{eq:F}
\end{equation}
and 
\begin{equation}
\underline{\rE}(r,t)<0 \quad \hbox{for any $t\in (0,\rT)$ and $a.e.~r\in(\varepsilon,\rR]$}.
\label{eq:sifno F}
\end{equation}
\label{prop:propauxiliar}
\end{prop}
\proof
Given an arbitrary long future horizon $\rT>0$ and each constant $\nu>1$ we introduce the function
$$
\underline{\rU}(|x|,t)=\Psi^{-1}\big (\nu (\rT-t +\rR- |x| \big )\big),\quad 0\le |x|<\rR,~0\le t<\rT.
$$
Clearly
$$
\left \{
\begin{array}{ll}
0<\underline{\rU}(|x|,t)<+\infty,& \quad \hbox{if $0\le |x|\le\rR,~0\le t<\rT$},\\ [.15cm]
\disp \lim_{t\nearrow\rT} \underline{\rU}(\rR,t)=+\infty. 
\end{array}
\right. 
$$
Straightforward computations on the function
$$
\underline{\rU}(|x|,t)=\Psi^{-1}(\zeta),
$$
for $\zeta =\nu\big (\rT-t +\rR-|x|\big )\in \big [\nu\big (\rT-t\big ),\nu (\rT-t)+\rR\big )\big ]$ and $0\le t<\rT$, shows that
$$
\begin{array}{ll}
-\Delta \underline{\rU}(|x|,t)+ g\big (\underline{\rU}(|x|,t)\big ) & \hspace*{-.3cm}=-\Delta
\Psi^{-1}(\zeta)+g\big (\Psi^{-1}(\zeta)\\ [.15cm]
& \hspace*{-.3cm}\le \big (1-\nu ^{2}\big )g\big (\Psi^{-1}(\zeta)\big )-\nu\dfrac{\rN-1}{|x|}\sqrt{2\rG\big (\Psi^{-1}(\zeta)\big)}<0
\end{array}
$$
(see \eqref{eq:tecnnivalityPsi}).
\par
\noindent
On the other hand, one has
$$
\left\{
\begin{array}[c]{l}
\underline{\rU}_{t}(\rR,t)=\nu \sqrt{2\rG\big (\Psi^{-1}\big (\nu_{1}\big (\rT-t)\big)\big )},\\[.2cm]
\disp \langle\nabla \underline{\rU} (\rR,t),\rn \rangle=\nu\sqrt{2\rG\big (\Psi^{-1}\big (\nu_{1}\big (\rT-t \big )\big)\big )},
\end{array}
\right.
$$
for $0\le t<\rT$. Thus
\begin{equation}
\underline{\rU}_{t}(\rR,t)+\dfrac{\partial}{\partial r}\underline{\rU}(\rR,t) =2\nu\sqrt{2\rG\big (\Psi^{-1}\big (\nu \big (\rT-t)\big)\big )},
\label{eq:technicalsubsolutionboundariesvalues}
\end{equation}
for $0\le t<\rT$. Under the assumptions \eqref{eq:dominationbalance},  the Proposition \ref{prop:wideregiondomination} implies
$$
2\nu \sqrt{2\rG\big ( \underline{\rU}(\rR,t)\big )}\le  \lambda f\big ( \underline{\rU}(\rR,t)  \big ),
$$ 
where the lower bound
$$
\Psi^{-1}\big (\nu (\rT+\rR\big )\big)\le \underline{\rU}(|x|,t),\quad 0\le |x|<\rR,~0\le t<\rT,
$$
is considered (see the comments of Remark \ref{rem:lambda0RT}). So that,  the inequality \eqref{eq:technicalsubsolutionboundariesvalues} leads to
$$
\underline{\rU}_{t}(\rR,t)+\dfrac{\partial}{\partial\rn}\underline{\rU}(\rR,t)\le   \lambda f\big (\underline{\rU}(\rR,t)\big ), \quad 0\le t<\rT.
$$ 
Finally, \eqref{eq:F} and \eqref{eq:sifno F} follow from the above arguments.$\fin$
\par
\medskip
Function $\underline{\rU}(r,t)$ looks like a subsolution but it is not in a strict sense since $\Delta
\underline{\rU}(\cdot, t)$ generates a measure, a Dirac mass at the origin (see Figure 2 below). The crucial fact to justify the comparison with the unique solution of $\rP(0)$ is that this measure is negative. A direct proof of the comparison can be given in term of  radially symmetric functions.
\begin{lemma} Assume  $\lambda >0$ satisfying the global domination of the Proposition \ref{prop:wideregiondomination}. Assume also  \eqref{eq:KellerOssermanG} and let $\rT=\Phi(u_{0})>0$. 
Let $u\in\cC\big ([0,\rT):\rH^{1}(0,\rR)\big )$ be the unique solution of \eqref{eq:uncontrolledradialproblem} and let $\underline{\rU}\in\cC\big ([0,\rT]:\cC^{\infty}(0, \rR)\big )$ defined by \eqref{eq:explosivesubsolution} for any $\nu>1$. Then
\begin{equation}
\underline{\rU}(r,t)\le u^{0}(r,t),\quad \hbox{for any $t\in [0,\rT[$ and any $r\in [0,\rR]$.}
\label{eq:radialcomparison}
\end{equation}
\label{lemma:comparison}
\end{lemma}
\proof By subtraction in the interior partial differential equations and multiplying by $\big (\underline{\rU}-u\big )_{+}$ we get, after an integration on the interval $(\varepsilon,\rR)$,
$$
\begin{array}{l}
\disp \dfrac{1}{2}\int_{\varepsilon}^{\rR}r^{\rN-1}\left (\dfrac{\partial }{\partial r}\big (\underline{\rU}-u^{0}\big )_{+}(r,t)\right )^{2}dr+\int_{\varepsilon}^{\rR}r^{\rN-1}\big (g(\underline{\rU}(r,t))-g(u^{0}(r,t))\big )\big (\underline{\rU}-u^{0}\big )_{+}(r,t)dr\\ [.3cm]
\disp =\int_{\varepsilon}^{\rR}r^{\rN-1}\underline{\rE}(r,t)dr
+\rR^{\rN-1}\left [\dfrac{\partial }{\partial r}\big (\underline{\rU}-u^{0}\big )(\rR,t)\right ]
\big (\underline{\rU}-u^{0}\big )_{+}(\rR,t)\\[.3cm]\quad  -
\varepsilon^{\rN-1}\left [\dfrac{\partial }{\partial r}\big (\underline{\rU}-u^{0}\big )(\varepsilon,t)\right ]
\big (\underline{\rU}-u^{0}\big )_{+}(\varepsilon,t).
\end{array}
$$
Since $\dfrac{\partial }{\partial r}\underline{\rU}(0,t)>0,~\dfrac{\partial }{\partial r}~u^{0}(0,t)=0$ and 
$$
\dfrac{\partial }{\partial r}~u^{0}(r,t)\ge 0,\quad r<\rR,~0<t<\rT,
$$
we deduce that there exists $\varepsilon>$ small enough for which
$$
\dfrac{\partial }{\partial r}\underline{\rU}(\varepsilon,t)>0\quad \hbox{for any $t\in [0,\rT]$}.
$$
So that 
$$
\liminf _{\varepsilon \searrow 0}\varepsilon^{\rN-1}\left [\dfrac{\partial }{\partial r}\big (\underline{\rU}-u^{0}\big )(\varepsilon,t)\right ]
\big (\underline{\rU}-u^{0}\big )_{+}(\varepsilon,t)\ge 0.
$$
Thus, we get
\vspace*{.2cm}
$$
\begin{array}{l}
\hspace*{-2cm}\disp \dfrac{1}{2}\left (\dfrac{\partial }{\partial t}\big (\underline{\rU}-u^{0}\big )_{+}(\rR,t)\right )^{2}+\dfrac{1}{2}\int_{0}^{\rR}r^{\rN-1}\left (\dfrac{\partial }{\partial r}\big (\underline{\rU}-u^{0}\big )_{+}(r,t)\right )^{2}dr\\ [.35cm]
\hspace*{2cm}\disp \le \lambda \int_{\varepsilon}^{\rR}\big (f(\underline{\rU}(\cdot,t))-f(u^{0}(\cdot,t))\big )\big (\underline{\rU}-u^{0}\big )_{+}(r,t)dr.
\end{array}
$$
for any $t\in [0,\rT)$. We choose $\rT=\Psi\big (u_{0}\big )$ for which
$$
\underline{\rU}(\rR,0)=\Psi^{-1}(\nu\rT)<\Psi^{-1}(\rT)=u_{0}.
$$
Being $f$ locally Lipschtiz, by applying  the Gronwall Lemma we have
\begin{equation}
\underline{\rU}(\rR,t)\le u^{0}(\rR,t),\quad 0<t<\rT.
\label{eq:inequalityinR}
\end{equation}
Finally, by repeating the integration by parts argument, now using \eqref{eq:inequalityinR}, we get
$$
\disp \dfrac{1}{2}\int_{\varepsilon}^{\rR}r^{\rN-1}\left (\dfrac{\partial }{\partial r}\big (\underline{\rU}-u^{0}\big )_{+}(\cdot,t)\right )^{2}dr\le 0,\quad 0<t<\rT
$$
and \eqref{eq:radialcomparison} follows.$\fin$
\bigskip
\par
\noindent
{\sc Proof of Theorem \ref{theo:boundaryblowup}}.\quad   With the notation of Proposition \ref{prop:propauxiliar}, Lemma \ref{lemma:comparison} implies 
$$
\Psi^{-1}\big (\nu(\rT-t +\rR- r\big )\big)\le u^{0}(r,t), \quad 0\le r<\rR,~0<t<\rT,
$$
with $\rT=\Psi \big (u_{0}\big )$. Since $\Psi (\infty)=0$ one deduces
$$
u^{0}(\rR,\rT)=+\infty.
$$ 
Therefore there exists a first boundary blow up time $\rT_{\infty}$ satisfying  $\rT_{\infty}\le \Psi \big (u_{0}\big )$ and 
$$
\dfrac{u^{0}(\rR,t)}{\Psi^{-1}\big (\nu\big (\rT_{\infty}-t)\big )\big)}\ge 1,\quad 0<t<\rT_{\infty}.
$$
Finally we argue as in \cite[Lemma 4.1]{AlDiRe}. We write the condition  \eqref{eq:nopotencias} as
\begin{equation}
\liminf_{s\rightarrow\infty} \frac{\Psi(\eta s)}{\Psi(s)}>1
\quad\hbox{for $ 0<\eta<1$}.
\label{eq:nopotenciasadecuada}
\end{equation}
(see Remark \ref{rem:nopotencias}). We may suppose that $0<(\nu-1)$ so small that
$$
\liminf_{s\rightarrow\infty} \frac{\Psi(\eta s)}{\Psi(s)}>\nu
\quad\hbox{for any $ 0<\eta<1$}.
$$
Thus if $0<\eta<1$ the assumption \eqref{eq:nopotenciasadecuada} implies  
$$
\Psi(\eta s)>\nu\Psi(s)\quad \hbox{for large $s$}.
$$
Therefore the monotonicity of $\Psi$ and
the choice $s=\Psi^{-1}\big (\rT_{\infty}-t\big )$ leads to
$$
\eta
\Psi^{-1}\big (\rT_{\infty}-t\big )<\Psi^{-1}\big (\nu\big (\rT_{\infty}-t \big )\big )
$$
and
$$
\dfrac{u^{0}(\rR,t)}{\Psi^{-1}(\rT_{\infty}-t)}\ge \dfrac{\Psi^{-1}\big (\nu(\rT_{\infty}-t)\big )}
{\Psi^{-1}\big (\rT_{\infty}-t)}>\eta,\quad 0<\rT_{\infty}-t\ll 1.
$$
Then
$$
\liminf_{t\nearrow\rT_{\infty}}\\
\dfrac{u^{0}(\rR,t)}{\Psi^{-1}\big (\rT_{\infty}-t\big )}>\eta. 
$$
Since $\eta<1$ is arbitrary one concludes the explosive boundary behaviour
$$
\liminf_{t\nearrow \rT_{\infty}}\dfrac{u^{0}(\rR,t)}{\Psi^{-1}\big (\rT_{\infty}-t\big )}\ge 1
$$
independently on $\nu>1.\fin$

\begin{rem}\rm 
Is is clear that the domination balance
$$
\lambda f(s)\ge \sqrt{2\rG(s)}\quad \hbox{for large $s$}
$$
implies 
$$
\Psi(s)=\int^{+\infty}_{s}\dfrac{ds}{\sqrt{2\rG(\tau)}}\le \dfrac{1}{\lambda} \int^{+\infty}_{s}\dfrac{d \tau}{f(\tau)} =\dfrac{1}{\lambda }\Phi(s)\quad \hbox{for large $s$}.
$$
Then
$$
\dfrac{\Psi \big (u^{0}(\rR,t)\big )}{\rT_{\infty}-t}\le \dfrac{1}{\lambda }\dfrac{\Phi \big (u^{0}(\rR,t)\big )}{\rT_{\infty}-t}\le 1,\quad \rT_{\infty}-t\ll 1
$$
and
$$
\limsup_{t\nearrow \rT_{\infty}}\dfrac{\Phi \big (u^{0}(\rR,t)\big )}{\rT_{\infty}-t}\le \lambda \quad \Rightarrow \quad \limsup_{t\nearrow \rT_{\infty}}\dfrac{\Psi \big (u^{0}(\rR,t)\big )}{\rT_{\infty}-t}\le 1.
$$ 
From suitable properties on $\Psi$ and $\Phi$, as  \eqref{eq:nopotencias}, one can deduce
$$
\liminf_{t\nearrow \rT_{\infty}}\dfrac{u^{0}(\rR,t)}{\Phi^{-1}\big (\lambda \big (\rT_{\infty}-t\big )\big )}\le 1 \quad \Rightarrow \quad \liminf_{t\nearrow \rT_{\infty}}\dfrac{u^{0}(\rR,t)}{\Psi^{-1}\big (\rT_{\infty}-t\big )}\le 1.
$$ 
\fineq
\label{rem:behaviouirbybalance}
\end{rem}
\par
\bigskip
The precise time growing rate of the trace of the solution $u(\rR,t)$ depends of the way in which the forcing term dominates over the absorption term. We will argue in the following by using  assumption \eqref{Hypo superlinear}.

\begin{rem}\rm As it was pointed out in Remark \ref{rem:behaviouirbybalance}, under \eqref{eq:dominationbalance} the condition \eqref{eq:KellerOssermanG} implies the superlinear condition \eqref{Hypo superlinear}. So that, under \eqref{eq:dominationbalance} the Remarks \ref{rem:potenciasg}, \ref{rem:exponential} and \ref{rem:exemplesg} and  (see also Remark \ref{rem:nopotencias})  provide some examples for which \eqref{Hypo superlinear} holds (see also Remark \ref{rem:exemplesf}).$\fin$
\end{rem}
\begin{figure}[htp]
	\begin{center}
		\includegraphics[width=14cm]{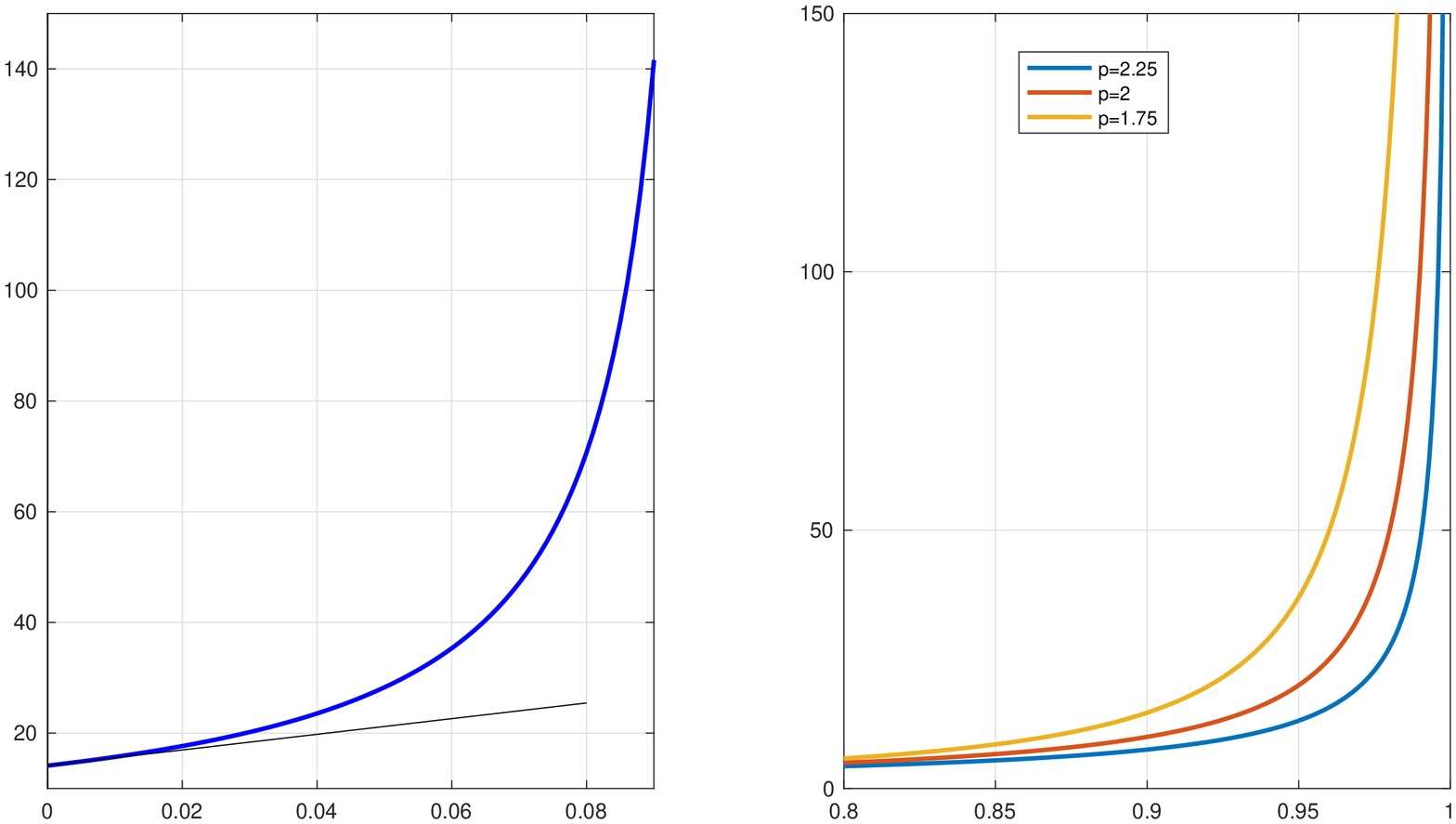}\\ 
		\hbox{\small Figure 2. Spatial profile of the subsolution $\underline{\rU}$ and time profiles of the solution at $r=\rR$ for some values of $p$.}
	\end{center}
\label{figure:profile}
\end{figure}
\par
We split the analysis in two different subsections.
\subsubsection{Time estimates for the strongly dominating  forcing over absorption case}
We assume, in this Subsection,  the condition (\ref{eq:fsuperwinG}). We have

\begin{theo}
Assume the hypothesis of Theorem \ref{theo:boundaryblowup}  as well as \eqref{eq:fsuperwinG}. 
Then the solution $u^{0}(r,t)$ of the problem \eqref{eq:uncontrolledradialproblem} verifies
$$
\liminf_{t\nearrow\rT_{\infty}}\dfrac{\Phi\big (u^{0}(\rR,t))}{\rT_{\infty}-t}\ge \lambda .
$$
More precisely, if we also assume 
\begin{equation}
\limsup_{s\rightarrow\infty} \frac{\Phi(\eta s)}{\Phi(s)} <1
\quad\hbox{for any } \eta>1. \label{eq:nopotenciasf}
\end{equation}
then we have the inequality
$$
\limsup_{t\nearrow\rT_{\infty}}\dfrac{u^{0}(\rR,t)}{\Phi^{-1}\big (\lambda (\rT_{\infty}-t)\big)}\le1. 
$$
\label{theo:blowupsuperdomination}
\end{theo}

\proof  We consider the interior equation 
$$
-\dfrac{1}{r^{\rN-1}}\dfrac{\partial}{\partial r}\left(  r^{\rN
-1}\dfrac{\partial u^{0}}{\partial r}(r,t)\right)  +g\big (u^{0}(r,t)\big )=0,\quad
r<\rR,~ 0<t<\rT_{\infty}.
$$
Multiplying by $r^{\rN-1}\dfrac{\partial u^{0}}{\partial r}$ one obtains
\begin{equation}
-\dfrac{1}{2}\dfrac{\partial}{\partial r}\left(  r^{\rN-1} \dfrac{\partial u^{0}}{\partial r}(r,t)\right)  ^{2}+r^{2(\rN-1)}
\dfrac{\partial}{\partial r}\rG\big (u^{0}(r,t)\big )=0,\quad
r<\rR,~ 0<t<\rT_{\infty}. 
\label{eq:radialmultiplied}
\end{equation}
Integrating from $0$ to $r<\rR$ we find
$$
0<\left(  r^{\rN-1}\dfrac{\partial u^{0}}{\partial r}(r,t)\right)
^{2}=2\int^{r}_{0}s^{2(\rN-1)}\dfrac{\partial}{\partial r}
\rG\big (u^{0}(s,t)\big )ds \le2 r^{2(\rN-1)}\rG\big (u^{0}(r,t)\big ).
$$
Thus
\begin{equation}
\dfrac{\partial u^{0}}{\partial r}(r,t)\le\sqrt{2\rG\big (u^{0}(r,t)\big )},\quad
r<\rR,~ 0<t<\rT_{\infty}.
\label{eq:boundu'}
\end{equation}
Therefore  the boundary condition leads to the ordinary differential inequality
\begin{equation}
\dfrac{\partial u^{0}}{\partial t}(\rR,t)+ \sqrt{2\rG\big (u^{0}(\rR,t)\big )}\ge \lambda f\big (u^{0}(\rR,t) \big ),\quad 0<t<\rT_{\infty}.
\label{eq:inequationdifferential}
\end{equation}
So that, applying assumption \eqref{eq:fsuperwinG} we deduce that for $0<\varepsilon<1$
$$
f\big (u^{0}(\rR,t)\big )\ge \varepsilon\sqrt{2\rG\big (u^{0}(\rR,t)\big )}\quad \hbox{for large values of $u^{0}(\rR,t)$}.
$$
and then
$$
\dfrac{\partial u^{0}}{\partial t}(\rR,t)\ge\lambda (1-\varepsilon)f\big (u^{0}(\rR,t)\big )\quad \hbox{for large values of $u^{0}(\rR,t)$},
$$
$i.e.$
$$
\dfrac{\dfrac{\partial u^{0}}{\partial t}(\rR,t)}{f\big (u^{0}\rR,t)\big )}\ge \lambda(1-\varepsilon), \quad\hbox{for large values of $u^{0}(\rR,t)$}.
$$
Integrating from $t$ to $\rT_{\infty}$ we get
\begin{equation}
\Phi\big (u^{0}(\rR,t)\big ) =\int^{\infty}_{u^{0}(\rR,t)}\dfrac
{ds}{f(s)}\ge \lambda (1-\varepsilon)(\rT_{\infty}-t),\quad0<\rT_{\infty}-t\ll1. \label{eq:casibehasup}
\end{equation}
Sending $\varepsilon\searrow0$ we obtain
$$
\liminf_{t\nearrow\rT_{\infty}}\dfrac{\Phi\big (u^{0}(\rR,t))}{\rT_{\infty}-t}\ge \lambda .
$$
We may make more precise this inequality as follows. From \eqref{eq:casibehasup} we have
\begin{equation}
u^{0}(\rR,t)\le\Phi^{-1}\big (\lambda (1-\varepsilon)(\rT_{\infty}-t)\big ),\quad0<\rT_{\infty}-t\ll1,
\label{eq:casicasibehasup}
\end{equation}
where the upper bound of the approach $\rT_{\infty}-t$ depends on $\varepsilon$ and it does not make to send $\varepsilon $ to 0.
Finally, given $\eta>1$ assumption \eqref{eq:nopotenciasf} implies
$$
\Phi(\eta s)<(1 -\varepsilon)\Phi(s)\quad\hbox{for large $s$}.
$$
Hence the monotonicity of $\Phi$ and the choice $s=\Phi^{-1}
\big (\lambda(\rT_{\infty}-t)\big)$ leads to
$$
\Phi^{-1}\big ((1-\varepsilon)\lambda (\rT_{\infty}-t)\big)< \eta\Phi
^{-1}\big (\lambda (\rT_{\infty}-t)\big),
$$
whence
$$
\dfrac{u^{0}(\rR,t)}{\Phi^{-1}\big (\lambda(\rT_{\infty}-t)\big)}\le
\dfrac{\Phi^{-1}\big ((1-\varepsilon)\lambda (\rT_{\infty}-t)\big )}{\Phi
^{-1}\big (\lambda (\rT_{\infty}-t)\big)}<\eta, \quad 0<\rT_{\infty}-t\ll1,
$$
and
$$
\limsup_{t\nearrow\rT_{\infty}(k)}\dfrac{u^{0}(\rR,t)}{\Phi
^{-1}\big (\lambda (\rT_{\infty}-t)\big)}\le\eta
$$
since $\eta>1$ is arbitrary the result holds.$\fin$

\begin{rem}\rm
If we consider the power-like choices $g_{m}(s)=s^{m}$ and $f_{p}(s)=s^{p}$
the assumptions of Theorem \ref{theo:blowupsuperdomination} hold whenever $p>\dfrac{m+1}{2}>1$. Then, since \eqref{eq:nopotenciasf} holds, one has
$$
\limsup_{t\nearrow \rT_{\infty}}u^{0}(\rR,t)\big (\rT_{\infty}-t\big )^{-\frac{1}{p-1}
}\le \left (\dfrac{1}{\lambda(p-1)}\right )^{\frac{1}{p-1}}
$$
We note that $u^{0}(\cdot,t)$ is integrable near $\rT_{\infty}$ if $p>2$. Once again, as in Remark \ref{rem:bang-bang}, this explains the non-uniqueness of the searched control.\fin
\label{rem:domination}
\end{rem}

\begin{rem}\rm
Under the technical assumption \eqref{eq:nopotencias} we deduce, as in
\cite[Theorem 1.1]{AlDiRe},
\begin{equation}
\lim_{|x|\nearrow\rR}\dfrac{u^{0}\big (|x|,\rT_{\infty}^{-}\big )}{\Psi
	^{-1}\big (\rR-|x|)}=1.\label{eq:profileballsolution}
\end{equation}
\fineqnum
\end{rem}
\begin{rem}\em Theorem \ref{theo:blowupsuperdomination}
also holds when we replace the assumptions \eqref{Hypo superlinear} and \eqref{eq:nopotenciasf} by 	\eqref{eq:monotocitysuperlinearf}. Indeed, from \eqref{eq:casicasibehasup} we obtain
$$
\dfrac{u^{0}(\rR,t)}{\Phi^{-1}\big (\lambda(\rT_{\infty}-t)\big)}\le
\dfrac{\Phi^{-1}\big ((1-\varepsilon)\lambda (\rT_{\infty}-t)\big )}{\Phi
^{-1}\big (\lambda (\rT_{\infty}-t)\big)}\le (1-\varepsilon)^{-\frac{1}{\alpha-1}}, \quad 0<\rT_{\infty}-t\ll1,
$$
by taking $\zeta=(1-\varepsilon)\lambda (\rT_{\infty}-t)$ and $\nu=(1-\varepsilon)^{-\frac{1}{\alpha-1}}>1$ in \eqref{eq:technicallityf} (see Remark \ref{rem:technicalityf}). Then 
$$
\limsup_{t\searrow \rT_{\infty}}\dfrac{u^{0}(\rR,t)}{\Phi^{-1}\big (\lambda(\rT_{\infty}-t)\big)}\le (1-\varepsilon)^{-\frac{1}{\alpha-1}}, \quad 0<\rT_{\infty}-t\ll1,
$$
and the result follows by letting $\varepsilon \searrow 0.\fin$
\label{rem:blowupsuperdominationmejora}
\end{rem}

\subsubsection{Time estimates for the weakly dominating forcing    over absorption case}

\begin{theo} Assume the hypothesis of Theorem \ref{theo:boundaryblowup}   as well as  \eqref{eq:fwinG-}. Then the solution of the problem \eqref{eq:uncontrolledradialproblem} behaves on the boundary as
$$
\liminf_{t\nearrow\rT_{\infty}}\dfrac{\Phi\big (u^{0}(\rR,t))}{\rT
_{\infty}-t}\ge\dfrac{\lambda \rL-1}{ \rL}.
$$
More precisely, under \eqref{eq:nopotenciasf} the inequality
$$
\limsup_{t\nearrow\rT_{\infty}}\dfrac{u^{0}(\rR,t)}{\Phi^{-1}\left(
\dfrac{\lambda\rL-1}{\rL}(\rT_{\infty}-t)\right)  }\le1,
$$
holds.
\label{theo:blowupdominationbalancedsup}
\end{theo}
\proof Arguing as in the proof of Theorem \ref{theo:blowupsuperdomination} we
obtain the ordinary differential inequality
$$
\dfrac{\partial u^{0}}{\partial t}(\rR,t)+ \sqrt{2\rG\big (u^{0}(\rR,t)\big )}\ge \lambda f\big (u^{0}(\rR,t) \big ),\quad 0<t<\rT_{\infty}
$$
(see \eqref{eq:inequationdifferential}).   So that, for $0<\varepsilon <\rL$  assumption \eqref{eq:fwinG-}, used in \eqref{eq:inequationdifferential}, leads to
$$
\dfrac{\partial u^{0}}{\partial t}(\rR,t)+\dfrac{1}{\rL-\varepsilon}f\big (u^{0}(\rR,t)\big )\ge \lambda f\big (u^{0}(\rR,t)\big )\quad
\hbox{for large values of $u^{0}(\rR,t)$},
$$
or equivalently
$$
\dfrac{\partial u^{0}}{\partial t}(\rR,t)\ge \dfrac{\lambda(\rL-\varepsilon)-1}{\rL-\varepsilon}f\big (u^{0}(\rR,t)\big ) \quad
\hbox{for large values of $u^{0}(\rR,t)$}. 
$$
Then
$$
\dfrac{\dfrac{\partial u^{0}}{\partial t}(\rR,t)}{f\big (u^{0}(\rR,t)\big )}\ge\dfrac
{\lambda(\rL-\varepsilon)-1}{\rL-\varepsilon}\quad
\hbox{for large values of $u^{0}(\rR,t)$},
$$
and 
\begin{equation}
\Phi\big (u^{0}(\rR,t)\big ) =\int^{\infty}_{u^{0}(\rR,t)}\dfrac
{ds}{f(s)}\ge\dfrac{\lambda(\rL-\varepsilon)-1}{\rL-\varepsilon}(\rT
_{\infty}-t),\quad0<\rT_{\infty}-t\ll1, \label{eq:casibehagenersup}
\end{equation}
by an integration from $t$ to $\rT_{\infty}$. Sending $\varepsilon\searrow0$ in \eqref{eq:casibehagenersup} we obtain
$$
\liminf_{t\nearrow\rT_{\infty}}\dfrac{\Phi\big (u^{0}(\rR,t))}{\rT
_{\infty}(\rR)-t}\ge \dfrac
{\lambda(\rL-\varepsilon)-1}{\rL-\varepsilon}
$$
and
$$
\liminf_{t\nearrow\rT_{\infty}}\dfrac{\Phi\big (u^{0}(\rR,t))}{\rT
	_{\infty}(\rR)-t}\ge \dfrac
{\lambda\rL-1}{\rL}.
$$
On the other hand, \eqref{eq:casibehagenersup} implies
\begin{equation}
u^{0}(\rR,t)\le\Phi^{-1}\left(\dfrac
{\lambda(\rL-\varepsilon)-1}{\rL-\varepsilon}(\rT_{\infty}-t)\right)  ,\quad0<\rT_{\infty}-t\ll1,
\label{eq:casicasibehasupL}
\end{equation}
where upper bound of the approach $\rT_{\infty}-t$ depends on $\varepsilon$ and it not has sense to send $\varepsilon $ to 0. Again we argue as in \cite[Lemma 4.1]{AlDiRe}. Given $\eta>1$ assumption \eqref{eq:nopotenciasf} implies
$$
\Phi(\eta s)<\dfrac{\rL\big (\lambda(\rL-\varepsilon)-1)\big )}{(\rL-\varepsilon)(\lambda \rL-1)}\Phi(s)\quad
\hbox{for large $s$}.
$$
Therefore the monotonicity of $\Phi$ and
the choice $s=\Phi^{-1}\left(\dfrac{\lambda \rL-1}{\rL}(
\rT_{\infty}-t)\right)  $ leads to
$$
\Phi^{-1}\left(\dfrac
{\lambda(\rL-\varepsilon)-1}{\rL-\varepsilon}
(\rT_{\infty}-t)\right)  < \eta\Phi^{-1}\left(\dfrac
{\lambda\rL-1}{\rL}(\rT_{\infty}-t)\right)
$$
whence
$$
\dfrac{u^{0}(\rR,t)}{\Phi^{-1}\left(\dfrac
{\lambda\rL-1}{\rL}\right)  }\le\dfrac{\Phi^{-1}\left(\lambda \dfrac{\rL -(\varepsilon+1)}{\rL-\varepsilon}(\rT_{\infty}-t)\right)
}{\Phi^{-1}\left( \dfrac
{\lambda\rL-1}{\rL}\right)
}<\eta, \quad0<\rT_{\infty}-t\ll1,
$$
and
$$
\limsup_{t\nearrow\rT_{\infty}}\dfrac{u^{0}(\rR,t)}{\Phi^{-1}\left(\dfrac
{\lambda\rL-1}{\rL}(\rT_{\infty}-t)\right)  }\le\eta.
$$
Since $\eta>1$ is arbitrary one concludes the result.$\fin$

\begin{rem}\rm
Once more, by \cite[Theorem 1.1]{AlDiRe} one deduces
$$
\lim_{|x| \nearrow \rR}\dfrac{u^{0}\big (|x|,\rT_{\infty}^{-}\big )}{\Psi
^{-1}\big (\rR-|x|\big )}=1,
$$
provided \eqref{eq:nopotencias}.$\fin$ 
\end{rem}
\begin{rem}\rm
For power-like choice  $g_{m}=s^{m}$,
with $m>1$, the assumption  \eqref{eq:fwinG-} becomes
$$
\liminf_{s\rightarrow\infty} f(s)s^{-\frac{m+1}{2}}=\rL\sqrt
{\dfrac{2}{m+1}},\quad \rL>0.
$$
In particular, for $f_{p}(s)=s^{p},~p>1,$ one has
$$
\Phi_{p}(s)=\dfrac{1}{(p-1)s^{p-1}},\quad s>0.
$$
Then, one obtains
$$
\limsup_{t\nearrow\rT_{\infty}}u^{0}(\rR,t)\big (\rT_{\infty
}-t\big )^{-\frac{1}{p-1}}\le \left(\dfrac{\rL}{(\lambda\rL-1)(p-1)}\right)  ^{\frac
{1}{p-1}},
$$
provided $p=\dfrac{m+1}{2}$ (see Theorem \ref{theo:invariance}). We note that $u^{0}(\cdot,t)$ is integrable near $\rT_{\infty}$ if $p>2.\fin$
\label{rem:power domination}
\end{rem}
\begin{rem}\em We may proceed  as in Remark 
\ref{rem:blowupsuperdominationmejora} in order to prove that 
Theorem \ref{theo:blowupdominationbalancedsup}
also holds when we replace the assumptions \eqref{Hypo superlinear} and \eqref{eq:nopotenciasf} by 	\eqref{eq:monotocitysuperlinearf}. Indeed, from \eqref{eq:casicasibehasupL} we obtain
$$
\dfrac{u^{0}(\rR,t)}{\Phi^{-1}\left(\dfrac
			{\lambda\rL-1}{\rL}(\rT_{\infty}-t)\right)  }\le\dfrac{\Phi^{-1}\left(\dfrac
{\lambda(\rL-\varepsilon)-1}{\rL-\varepsilon}(\rT_{\infty}-t)\right)}
{\Phi^{-1}\left(\dfrac
	{\lambda\rL-1}{\rL}(\rT_{\infty}-t)\right) }\le
\left (\dfrac
{\lambda\rL-1}{\rL}\dfrac
	{\rL-\varepsilon}{\lambda(\rL-\varepsilon)-1}\right )^{\frac{1}{\alpha -1}},
$$	
for $0<\rT_{\infty}-t\ll1$, by taking $\zeta=\dfrac
{\lambda(\rL-\varepsilon)-1}{\rL-\varepsilon}(\rT_{\infty}-t)$ and $\nu=\left (\dfrac
{\lambda\rL-1}{\rL}\dfrac
{\rL-\varepsilon}{\lambda(\rL-\varepsilon)-1}\right )^{\frac{1}{\alpha -1}}$ in \eqref{eq:technicallityf} (see Remark \ref{rem:technicalityf}). Since the function
$$
\varepsilon \mapsto \dfrac{\lambda(\rL-\varepsilon)-1}
{\rL-\varepsilon}
$$
is decreasing, we have $\nu>1$. Then 
$$
\limsup_{t\searrow \rT_{\infty}}\dfrac{u^{0}(\rR,t)}{\Phi^{-1}\left(\dfrac
	{\lambda\rL-1}{\rL}(\rT_{\infty}-t)\right)  }\le
\left (\dfrac
{\lambda\rL-1}{\rL}\dfrac
{\rL-\varepsilon}{\lambda(\rL-\varepsilon)-1}\right )^{\frac{1}{\alpha -1}},
$$	
and the result follows by letting $\varepsilon \searrow 0.\fin$
\label{rem:blowupdominationbalancedsupmejora}
\end{rem}
\par

The study of the behaviour at the finite blow up time is completed now as it is collected in the following result.

\begin{theo}[Behaviour at the finite blow up time]Suppose 
\begin{equation}
\liminf_{\tau\rightarrow\infty}\dfrac{f(\tau)}{\sqrt{2\rG(\tau)}}
=\ell\label{eq:fwinG+}>1
\end{equation}
and the assumptions of Theorem \ref{theo:boundaryblowup} .
Then the boundary behaviour of the solution, $u^{0}(r,t)$ of \eqref{eq:uncontrolledradialproblem} verifies 
$$
\liminf_{t\nearrow\rT_{\infty}}\dfrac{u^{0}(\rR,t)}{\Phi^{-1}\left(
\dfrac{\lambda \ell-1}{\ell}(\rT_{\infty}-t)\right)  }\ge 1.
$$
If \eqref{eq:fwinG-} also holds one has 
\begin{equation}
\limsup_{t\nearrow\rT_{\infty}}\dfrac{u^{0}(\rR,t)}{\Phi^{-1}\left(
	\dfrac{\lambda \rL-1}{\rL}(\rT_{\infty}-t)\right)  }\le1\le
\liminf_{t\nearrow\rT_{\infty}}\dfrac{u^{0}(\rR,t)}{\Phi^{-1}\left( 
	\dfrac{\lambda \ell-1}{\ell}(\rT_{\infty}-t)\right)  }.
\label{eq:estiboundary}
\end{equation}
In particular, the property
\begin{equation}
\lim_{\tau\rightarrow\infty}\dfrac{f(\tau)}{\sqrt{2\rG(\tau)}}=\ell >1
\label{eq:datosuniinfty}
\end{equation}
implies
$$
\lim_{t\nearrow\rT_{\infty}}\dfrac{u^{0}(\rR,t)}{\Phi^{-1}\left(\dfrac
{\lambda \ell-1}{\ell}(\rT_{\infty}-t)\right)  }=1.
$$ 
\label{theo:blowupdominationbalanced}
\end{theo}
\proof First of all we note that from Theorem \ref{theo:blowupdominationbalancedsup} one satisfies
$$
\limsup_{t\nearrow\rT_{\infty}}\dfrac{u^{0}(\rR,t)}{\Phi^{-1}\left( 
\dfrac{\lambda \rL-1}{\rL}(\rT_{\infty}-t)\right)  }\le1.
$$
Once more, we recall that under \eqref{eq:fwinG+} the assumption \eqref{eq:KellerOssermanG} implies \eqref{Hypo superlinear}.
\par
\noindent In order to complete \eqref{eq:estiboundary}
we use the inequality
$$
\dfrac{\partial u^{0}}{\partial r}(r,t)\le\sqrt{2\rG\big (u^{0}(r,t)\big )},\quad r<\rR,~0<t<\rT_{\infty}
$$
(see \eqref{eq:boundu'}) on the interior equation
$$
\dfrac{\partial^{2}u^{0}}{\partial r^{2}}(r,t)+\dfrac{\rN-1}{r}
\dfrac{\partial u^{0}}{\partial r}(r,t)=g\big (u^{0}(r,t)\big ),\quad r<\rR,~0<t<\rT_{\infty}.
$$
Then
$$
\dfrac{\partial^{2}}{\partial r^{2}}u^{0}(r,t)+\dfrac{\rN-1}{r}
\sqrt{2\rG\big (u^{0}(r,t)\big )}\ge g\big (u^{0}(r,t)\big ),\quad
0<r<\rR,~0<t<\rT_{\infty}.
$$
On the other hand, by using \eqref{eq:potenciamayorqueuno} of Lemma
\ref{Lemma:KellerOssermanAGfraccion} we deduce for each $\varepsilon>0$
$$
\dfrac{\partial^{2}}{\partial r^{2}}u^{0}(r,t)+\varepsilon \dfrac{\rN-1}{r}g\big (u^{0}(r,t)\big )\ge g\big (u^{0}(r,t)\big )
$$
whenever
$$
0<\big (\rT_{\infty}-t\big )+\big (\rR-r\big )\le\delta\ll1.
$$
In particular, as $\delta <\dfrac{\rR}{2}$ we have $r>\dfrac{\rR}{2}$ and 
$$
\dfrac{\partial^{2}}{\partial r^{2}}u^{0}(r,t)\ge
\left(  1-2\varepsilon \dfrac{\rN-1}{\rR}\right) g\big (u^{0}(r,t)\big ).
$$
Multiplying by $\dfrac{\partial}{\partial r}u^{0}(r,t)>0$ in these region we deduce
$$
\dfrac{1}{2}\dfrac{\partial}{\partial r}\left(  \dfrac{\partial}{\partial
r}u^{0}(r,t)\right)  ^{2}\ge\left(  1-2\varepsilon \dfrac{\rN-1}{\rR}\right)  \dfrac{\partial}{\partial r}\rG\big (u^{0}(r,t)\big ).
$$
Next, an integration from $\rR-\delta$ to $\rR$ leads to
$$
\dfrac{\partial}{\partial r}u^{0}(\rR,t)\ge\sqrt{2\left( 1-2\varepsilon
\dfrac{\rN-1}{\rR}\right)  \big (\rG\big (u^{0}(\rR,t)\big )-\rG\big (u^{0}(\rR-\delta
,t)\big )\big )}.
$$
Since $\rT_{\infty}$ is the first time in which  $u(\rR,\cdot)$ becomes infinity, we may suppose $t$ so close to $\rT_{\infty}$ as
$$
\rG\big (u^{0}(\rR,t)\ge \dfrac{1}{\varepsilon}\rG\big (\rU_{\infty
}^{\rB_{\rR}}(\rR-\delta)\big ).
$$
Hence
$$
\rG\big (u^{0}(\rR,t)\ge \dfrac{1}{\varepsilon}\rG\big (\rU_{\infty
}^{\rB{\rR}}(\rR-\delta)\big )>\dfrac{1}{\varepsilon}\rG\big (u^{0}(\rR-\delta,t)\big ),
$$
provided $t$ near $\rT_{\infty}$ and $\delta$ small (here by $\rU_{\infty}^{\bB_{\rR}}$ we are denoting the relative radial symmetric large solution of  \eqref{eq:largesolution}). It implies
$$
\dfrac{\partial}{\partial r}u^{0}(\rR,t)\ge\sqrt{2\left(  1-2\varepsilon \dfrac{\rN-1}{\rR}\right)  (1-\varepsilon)\rG \big (u^{0}(\rR,t)\big )},\quad 0<\rT_{\infty}-t\ll 1,
$$
consequently from the boundary condition we get to  the ordinary differential inequality
$$
\dfrac{\partial u^{0}}{\partial t}(\rR,t)+\sqrt{2\left( 1-2\varepsilon \dfrac{\rN-1}{\rR}\right ) (1-\varepsilon)
\rG\big (u^{0}(\rR,t)\big )}\le \lambda f\big (u^{0}(\rR,t) \big ).
$$
Since assumption \eqref{eq:fwinG+} implies 
$$
\ell -\varepsilon\le \dfrac{f\big (u^{0}(\rR,t) \big )}{\sqrt{2\rG\big (u^{0}(\rR,t)\big )}}\le \ell +\varepsilon\quad \hbox{for large values of $u^{0}(\rR,t)$},
$$
we deduce
$$
\dfrac{\partial u^{0}}{\partial t}(\rR,t)\le \left  (\dfrac{\lambda(\ell +\varepsilon) -\sqrt{\left( 1-2\varepsilon \dfrac{\rN-1}{\rR}\right ) (1-\varepsilon)}}{\ell +\varepsilon}\right )f\big (u^{0}(\rR,t) \big ) \quad \hbox{for large values of $u^{0}(\rR,t)$}
$$
and
$$
\dfrac{\dfrac{\partial u^{0}}{\partial t}(\rR,t)}{f\big (u^{0}(\rR
,t)\big )}\le\dfrac{\lambda \ell +\rE(\varepsilon)}{\ell +\varepsilon}\quad \hbox{for large values of $u^{0}(\rR,t)$}
$$
where the function  $\rE(\varepsilon)\doteq \lambda \varepsilon -\sqrt{\left( 1-2\varepsilon \dfrac{\rN-1}{\rR}\right ) (1-\varepsilon)}$ satisfies $\disp \lim_{\varepsilon \rightarrow 0}\rE(\varepsilon)=-1$. So that, an integration from $t$ to
$\rT_{\infty}$ leads to
$$
\Phi\big (u^{0}(\rR,t)\big )=\int^{\infty}_{u^{0}(\rR,t)}\dfrac
{ds}{f(s)}\le \dfrac{\lambda \ell +\rE(\varepsilon)}{\ell +\varepsilon}(\rT_{\infty}-t)
$$
whence 
$$
\lim_{t\nearrow \rT_{\infty}}\dfrac{\Phi\big (u^{0}(\rR,t)\big )}{\rT_{\infty}-t}\le  \dfrac{\lambda \ell +\rH(\varepsilon)}{\ell +\varepsilon}
$$
and 
$$
\lim_{t\nearrow \rT_{\infty}}\dfrac{\Phi\big (u^{0}(\rR,t)\big )}{\rT_{\infty}-t}\le \dfrac{\lambda \ell -1}{\ell}
$$
by letting $\varepsilon \rightarrow 0$. Once more, we note that in the inequality
\begin{equation}
u^{0}(\rR,t)\ge\Phi^{-1}\left (\dfrac{\lambda \ell +\rE(\varepsilon)}{\ell +\varepsilon}(\rT_{\infty}-t) \right )  ,\quad 0<\rT_{\infty}-t\ll 1
\label{eq:casicasilast}
\end{equation}
the upper bound of the approach $\rT_{\infty}-t$ depends on $\varepsilon$ and it does not make sense to send $\varepsilon $ to 0.
However we may argue by using the assumption \eqref{eq:nopotenciasf}. Indeed, reasoning as in  Remark \ref{rem:nopotencias} we may write \eqref{eq:nopotenciasf} as
\begin{equation}
\liminf_{s\rightarrow\infty} \frac{\Phi(\eta s)}{\Phi(s)}>1\quad \hbox{for any $ 0<\eta<1$}.
\label{eq:nopotenciasfadecuada}
\end{equation}
Since $\dfrac{\ell \big (\lambda \ell+\rE(\varepsilon \big ))}{(\ell +\varepsilon)(\lambda \ell -1)}>1$, for $\varepsilon$ small,  given $0<\eta <1$,  \eqref{eq:nopotenciasfadecuada}  implies
$$
\dfrac{\ell \big (\lambda \ell+\rE(\varepsilon \big ))}{(\ell +\varepsilon)(\lambda \ell -1)}\Phi(s)<\Phi(\eta s)\quad \hbox{for large $s$}.
$$
Consequently the monotonicity of $\Phi$ and the choice $s=\Phi^{-1}\left (\dfrac{\lambda \ell -1}{\ell }(\rT_{\infty}-t) \right) 
$ leads to
$$
\eta\Phi^{-1}\left (\dfrac{\lambda \ell -1}{\ell}(\rT_{\infty}-t) \right ) <\Phi^{-1}\left (\dfrac{\lambda \ell +\rE(\varepsilon)}{\ell +\varepsilon}(\rT_{\infty}-t) \right ),
$$
whence
$$
\dfrac{u^{0}(\rR,t)}{\Phi^{-1}\left (\dfrac{\lambda \ell -1}{\ell }(\rT_{\infty}-t) \right)  }\ge\dfrac{\Phi^{-1}\left (\dfrac{\lambda \ell +\rE(\varepsilon)}{\ell +\varepsilon}(\rT_{\infty}-t) \right ) }{\Phi^{-1}\left (\dfrac{\lambda \ell -1}{\ell}(\rT_{\infty}-t) \right )  }>\eta,
$$
for $\rT_{\infty}-t$ small, and we have
$$
\liminf_{t\nearrow\rT_{\infty}}\dfrac{u^{0}(\rR,t)}{\Phi^{-1}\left (\dfrac{\lambda \ell -1}{\ell }(\rT_{\infty}-t) \right)  }\ge\eta
$$
since $\eta<1$ is arbitrary one concludes the result.$\fin$
\begin{rem}\rm For power-like choice  $g_{m}=s^{m}$,
with $m>1$, the assumption \eqref{eq:datosuniinfty} becomes
$$
\lim_{s\rightarrow\infty} f(s)s^{-\frac{m+1}{2}}=\ell\sqrt
{\dfrac{2}{m+1}},~\ell>1.
$$
In particular, for $f_{p}(s)=s^{p},~p>1$ one has
$$
\Phi_{p}(s)=\dfrac{1}{(p-1)s^{p-1}},\quad s>0.
$$
Then, one obtains
$$
\lim_{t\nearrow\rT_{\infty}}u^{0}(\rR,t)\big (\rT_{\infty
}-t\big )^{-\frac{1}{p-1}}=\left(\dfrac{\ell}{(\lambda \ell-1)(p-1)}\right)  ^{\frac
{1}{p-1}},
$$
provided $p=\dfrac{m+1}{2}$.
We note that $u^{0}(\cdot,t)$ is integrable near $\rT_{\infty}$ if $p>2$.
\end{rem}
\begin{rem}\rm As in Remark \ref{rem:nopotencias} one may extend the arguments  to a kind of borderline case given by
$$
\limsup_{s\rightarrow\infty} \frac{\Phi(\eta_{0} s)}{\Phi(s)}=1
\quad\hbox{for some $\eta_{0}>1$}, 
$$
but by simplicity we will do not consider that in this paper.$\fin$
\label{rem:nopotenciasfborder}
\end{rem}
\begin{rem}\rm Once more, we may prove Theorem \ref{theo:blowupdominationbalanced} when we replace the assumptions \eqref{Hypo superlinear} and \eqref{eq:nopotenciasf} by 	\eqref{eq:monotocitysuperlinearf}. Indeed, from \eqref{eq:casicasilast} we obtain
$$
\dfrac{u^{0}(\rR,t)}{\Phi^{-1}\left (\dfrac{\lambda \ell -1}{\ell }(\rT_{\infty}-t) \right)  }\ge\dfrac{\Phi^{-1}\left (\dfrac{\lambda \ell +\rE(\varepsilon)}{\ell +\varepsilon}(\rT_{\infty}-t) \right ) }{\Phi^{-1}\left (\dfrac{\lambda \ell -1}{\ell}(\rT_{\infty}-t) \right )  }\ge \left (\dfrac{\lambda \ell -1}{\ell }\dfrac{\ell +\varepsilon}{\lambda \ell +\rE(\varepsilon)}\right )^{\frac{1}{\alpha -1}},
$$
for $\rT_{\infty}-t$ small, by taking $\zeta=\dfrac{\lambda \ell +\rE(\varepsilon)}{\ell +\varepsilon}(\rT_{\infty}-t)$ and $\nu=\left (\dfrac{\lambda \ell -1}{\ell }\dfrac{\ell +\varepsilon}{\lambda \ell +\rE(\varepsilon)}\right )^{\frac{1}{\alpha -1}}$ in \eqref{eq:technicallityf2} (see Remark \ref{rem:technicalityf}). We note that $\nu<1$ as it was pointed out in the above proof of Theorem~\ref{theo:blowupdominationbalanced}. Then
$$
\liminf_{t\searrow \rT_{\infty}}\dfrac{u^{0}(\rR,t)}{\Phi^{-1}\left (\dfrac{\lambda \ell -1}{\ell }(\rT_{\infty}-t) \right)  }\ge\left (\dfrac{\lambda \ell -1}{\ell }\dfrac{\ell +\varepsilon}{\lambda \ell +\rE(\varepsilon)}\right )^{\frac{1}{\alpha -1}},
$$
and the result follows by letting $\varepsilon \searrow 0$.
\par
\noindent
We emphasize that this reasoning, as well as the ones of Remarks \ref{rem:blowupsuperdominationmejora} and \ref{rem:blowupdominationbalancedsupmejora}, includes the borderline case announced in Remark \ref{rem:nopotenciasfborder}.$\fin$
\label{rem:blowupdominationbalancedmejora}
\end{rem}
\par
 We consider now the limiting (weak domination) case $p=(m+1)/2$, appearing in the power-like problem governed by the choices $g_{m}(r)=r^{m},~m>0,$ and $f_{p}(r)=r^{p},~p>0$, of the domination assumption. To better illustrate the behaviour near the finite blow up time $\rT_{\infty}$, we consider the case of the self-similar solution corresponding to the spatial domain given by the hyperplane $\RR^{\rN-1}\times \RR_{+}$. This special case provides some intrinsic information in this limit case. We have
\begin{theo} Assume $2p=m+1$. Then, the system
$$
\left\{
\begin{array}{ll}
-\Delta u+u^{m}=0 & \quad\hbox{in }\big (\RR^{\rN-1}\times \RR_{+}\big )\times(0,\infty),\\[.15cm]%
\dfrac{\partial u}{\partial t}+\dfrac{\partial u}{\partial\bn }=u^{p} &
\quad\hbox{on }\big (\RR^{\rN-1}\times \{0\}\big )\times (0,\infty),
\end{array}
\right.  
$$
is invariant by the change of variables $v(x,t)=\mu ^{q}u(\mu x,\mu t),~\mu >0$. In this case, any self-similar solution must be of the form
\begin{equation}
u(x,t)=\dfrac{1}{t^{\frac{1}{p-1}}}\rH\left (\dfrac{x}{t}\right ),\quad x\in\RR^{\rN-1}\times \RR_{+},~0<t,
\label{eq:selfsimilarity}
\end{equation}
for some $p\neq 1$ and some function $\rH :\RR^{\rN-1}\times \RR_{+}\rightarrow \RR$, called as the similarity profile of the self-similar solution, satisfying the boundary value problem
\begin{equation}
\left \{
\begin{array}{l}
-\Delta \rH (\eta)+\big (\rH(\eta)\big )^{m}=0,\quad \eta\in \RR^{\rN-1}\times \RR_{+},\\ [.2cm]
\disp \sum_{i=1}^{\rN-1}\eta_{i}\rD_{i}\rH(\eta)-\rD_{\rN}\rH(\eta)=\dfrac{1}{p-1}\rH(\eta)+\big (\rH(\eta)\big )^{p},\quad \eta\in \RR^{\rN-1}\times \{0\}.
\end{array}
\right .
\label{eq:similaityprofile}
\end{equation}
In particular, the function with a regional blowing-up set
\begin{equation}
\rH(\eta)=\left (\dfrac{2(m+1)}{(m-1)^{2}}\right )^{\frac{1}{m-1}}\left [\eta_{\rN}-\dfrac{1}{2}\big (\sqrt{2(m+1)}-2\big )
\right ]_{+}^{-\frac{2}{m-1}},\quad \eta\in\RR^{\rN-1}\times \RR_{+}
\label{eq:similarityprofile}
\end{equation}
is the similarity profile of the self-similar solution
\begin{equation}
u(x,t)=\left (\dfrac{2(m+1)}{(m-1)^{2}}\right )^{\frac{1}{m-1}}\left [x_{\rN}-t\dfrac{\sqrt{2(m+1)}-2}{2}
\right ]_{+}^{-\frac{2}{m-1}},
\label{eq:selfsimilaritydistance}
\end{equation}
for $x=(x',x_{\rN})\in\RR^{\rN-1}\times \RR_{+},~0<t$, 
where expression \eqref{eq:selfsimilaritydistance} is uniform on $x'\in\RR^{\rN-1}$.
\label{theo:invariance}
\end{theo}
\proof  The function $v(x,t)=\mu ^{q}u(\mu x,\mu t),~\mu >0$, verifies
$$
\left \{
\begin{array}{l}
-\Delta v+g(v)=\mu ^{2-q(m-1)}v^{m},\\ [.2cm]
\dfrac{\partial v}{\partial t}-\dfrac{\partial v}{\partial x_{\rN}}=\mu ^{1-q(p-1)}v^{p}.
\end{array}
\right.
$$
So that, the invariance of the equations follows if 
$$
2-q(m-1)=0\quad\hbox{and} \quad 1-q(p-1)=0,
$$
whence $2p=m+1$. The self-similarity equality is 
$$
u(x,t)=\mu ^{\frac{2}{p-1}}u(\mu x,\mu t).
$$
Then, if we take $\mu =1$, after a derivation  with respect to $\mu$, one concludes
$$
\pe{x}{\nabla u(x,t)}+tu_{t}(x,t)+\dfrac{1}{p-1}u(x,t)=0.
$$
Then, by means of the classical characteristics method, we obtain that any self-similar solutions is represented by
$$
u\big (xe^{s},te^{s}\big )=u(x,t)e^{-\frac{1}{p-1}s},\quad s>0,
$$
whence
$$
u(x,t)=\dfrac{1}{t^{\frac{1}{p-1}}}\rH\left(\frac{x}{t}\right ),\quad (x,t)\in\big (\RR^{\rN -1}\times \RR_{+}\big )\times (0,\infty)
$$ 
for a profile function $\rH:\RR^{\rN-1}\times \RR_{+}\rightarrow \RR$ verifying
$$
\left \{
\begin{array}{l}
-\Delta \rH (\eta)+\big (\rH(\eta)\big )^{m}=0,\quad \eta\in \RR^{\rN-1}\times \RR_{+},\\ [.2cm]
\disp \sum_{i=1}^{\rN-1}\eta_{i}\rD_{i}\rH(\eta)-\rD_{\rN}\rH(\eta)=\dfrac{1}{p-1}\rH(\eta)+\big (\rH(\eta)\big )^{p},\quad \eta\in \RR^{\rN-1}\times \{0\}.
\end{array}
\right .
$$
Let us prove that we can take, as a possible solution, the profile function given by
$$
\rH(\eta)=k\big (\eta_{\rN}-\rC\big )^{\alpha},\quad \eta\in \RR^{\rN-1}\times \RR_{+}
$$
for some suitable positive constants $k,~\alpha$ and $\rC$. Indeed, clearly
$$
-\Delta \rH (\eta)+\big (\rH(\eta)\big )^{m}=0\quad \Leftrightarrow \quad 
\rH(\eta)=k_{m}\big (\eta_{\rN}-\rC\big )^{-\frac{2}{m-1}},\quad \eta_{\rN}>0
$$
and $k_{m} =\left (\dfrac{2(m+1)}{(m-1)^{2}}\right )^{\frac{1}{m-1}}$. Moreover, the boundary condition becomes
$$
k_{m}^{p}(-\rC)^{-\frac{2p}{m-1}}+\dfrac{k_{m}}{p-1}(-\rC)^{-\frac{2}{m-1}}-
\dfrac{2k_{m}}{m-1}(-\rC)^{-\frac{2}{m-1}-1}=0,
$$
and thus we must require
$$
k_{m}^{p}(-\rC)^{-\frac{2(p-1)}{m-1}}+\dfrac{k_{m}}{p-1}-
\dfrac{2k_{m}}{m-1}(-\rC)^{-1}=0.
$$
Since $2p=m+1$ we know that $\dfrac{2(p-1)}{m-1}=\dfrac{m+1-2}{m-1}=1$ and we have
$$
\left (k_{m}^{p-1}-\dfrac{2}{m-1}\right )(-\rC)^{-1}=-\dfrac{1}{p-1},
$$
whence
$$
\rC=(p-1)\left (k_{m}^{p-1}-\dfrac{2}{m-1}\right )=\dfrac{1}{2}\big (\sqrt{2(m+1)}-2\big ).
$$
\fineqnum
\vspace*{-.5cm}
\begin{figure}[htp]
	\begin{center}
		\includegraphics[width=14cm]{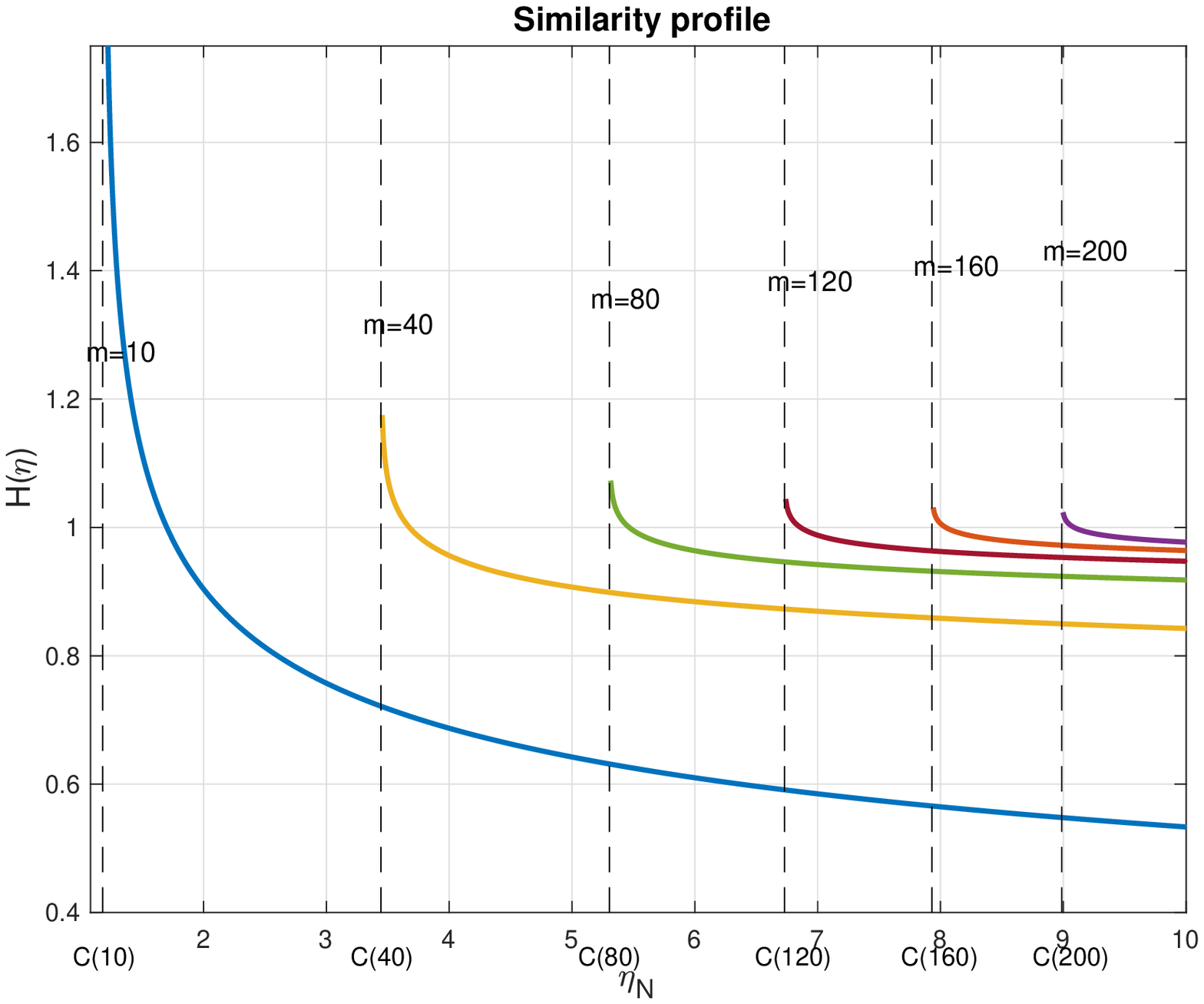}
	\end{center}
	\label{figure:similarityprofile}
\end{figure}
\begin{rem}[Blow up finite time property]\rm Since $2p=m+1,~p\neq 1$, the self-similar solution defined by \eqref{eq:selfsimilaritydistance} can be represented as
\begin{equation}
u(x,t)=\left (\dfrac{\sqrt{p}}{(p-1)(\sqrt{p}-1)}\right )^{\frac{1}{p-1}}
\left [\dfrac{x_{\rN}}{\sqrt{p}-1}-t
\right ]_{+}^{-\frac{1}{p-1}},
\label{eq:selfsimilarityblowup}
\end{equation}
for $\quad x=(x',x_{\rN})\in \RR^{\rN-1}\times \RR_{+},~t>0$.
Then for each $x_{\rN}>0$ the solution at this point blows up at the finite time 
\begin{equation}
\rT_{\infty}(x_{\rN})\doteq \dfrac{x_{\rN}}{\sqrt{p}-1},
\label{eq:blowupfinitetime}
\end{equation}
i.e.,
$$
\left \{
\begin{array}{ll}
0<u(x,t)<+\infty & \hbox{if $t<\rT_{\infty}(x_{\rN})$},\\ [.15cm]
u(x,t)=+\infty  & \hbox{if $t\ge \rT_{\infty}(x_{\rN})$}.
\end{array}
\right .
$$
Notice that, again \eqref{eq:selfsimilarityblowup} is uniform on $x'\in\RR^{\rN-1}.\fin$
\label{rem:blowupsimilarity}
\end{rem}
\begin{rem}[Large solution at the blow up time]\rm Consider now the spatial domain $\RR^{\rN-1}\times (\rR,+\infty),~\rR>0.$ From \eqref{eq:selfsimilaritydistance} we obtain the representation 
\begin{equation}
u\big (x,\rT_{\infty}(\rR)\big )=\left (\dfrac{2(m+1)}{(m-1)^{2}}\right )^{\frac{1}{m-1}}
\left [x_{\rN}-\rR
\right ]_{+}^{-\frac{2}{m-1}},
\label{eq:selfsimilaritylargesolution}
\end{equation}
for $\quad x=(x',x_{\rN})\in \RR^{\rN-1}\times (\rR,+\infty)$, where
$$
\rT_{\infty}(\rR)\doteq \dfrac{\rR}{\sqrt{p}-1}=\dfrac{2\rR}{\sqrt{2(m+1)}-2}.
$$
Thus
$$
\left \{
\begin{array}{ll}
0<u\big (x,\rT_{\infty}(\rR)\big )<+\infty & \hbox{if $x_{\rN}>\rR$},\\ [.15cm]
u(x,t)=+\infty  & \hbox{if $x_{\rN}=\rR$}.
\end{array}
\right .
$$
Once more, \eqref{eq:selfsimilaritylargesolution} is uniform on $x'\in\RR^{\rN-1}$. Notice that we have
that 
$$
\dfrac{\partial u}{\partial t}(\rR,t)+\dfrac{\partial u}{\partial \bn}(\rR,t)-u^{p}(\rR,t)=\gamma(t),\quad 0<t<\rT_{\infty}(\rR),
$$
where
$$
\gamma (t)=\left (\dfrac{\sqrt{p}}{(p-1)(\sqrt{p}-1)}\right )^{\frac{1}{p-1}}
\dfrac{\sqrt{p}-(p+1)}{(p-1)(\sqrt{p}-1)}\left [\rT_{\infty}(\rR)-t
\right ]_{+}^{-\frac{p}{p-1}}.
$$
\fineq
\label{rem:largesolutionsimilarity}
\end{rem}
\begin{rem}\rm Notice that the above boundary condition can be equivalently formulated as an "oblique boundary" condition with a nonlinear forcing term which, as far as we know, was not previously treated in the literature. Some techniques and references on nonlinear oblique problems can be found in~\cite{DDO}.$\fin$
\end{rem}

\subsection{Proof of the controlled explosions for problem $\rP(\alpha)$}

Now we return to the treatment of the problem $\rP(\alpha)$.
\medskip

\begin{lemma} 
\label{lemma:compariso P_0}
Let $u^0(r,t)$ be the unique solution of problems $\rP(0)$ with blow-up time on the boundary $\rT_{\infty}$ and let a given $\varepsilon>0$.  Let $\rV^{\alpha }(t)$ be the solution of problem $\rP_{\rF(\cdot,\alpha )}$, given in Theorem 1,  corresponding to the control $\alpha (t)$ associated to the time $\rT_{\infty}$ and initial datum  $\rV^{\alpha }(0) \geq u_{0}$ . Then we have the comparison
$$
0\leq u^{\alpha }(t,x)\leq \rV^{\alpha }(t)\text{ for any }t\in [
0,\rT_{\infty})\text{ and any }x\in \rB_{\rR}. 
$$
\end{lemma}
\noindent Proof.  We know that $\rV^{\alpha }(t)\geq 0$ for almost every $t>0$. Then the function $
w^{\alpha }(t,x)=\rV^{\alpha }(t)$ satisfies that 
$$
\left\{ 
\begin{array}{ll}
-\Delta w^{\alpha }+\mathcal{G}(w^{\alpha },\alpha )\geq 0 & \text{for }(x,t)\in  \bB_{\rR}\times  (0,\rT_{\infty}), \\ [.2cm]
\dfrac{\partial w^{\alpha }}{\partial t}+\dfrac{\partial w^{\alpha }}{\partial
\bn}= \lambda f_{\rM_{\varepsilon}}(u)+\alpha (t) & \text{for }(x,t)\in 
\partial \bB_{\rR}\times (0,\rT_{\infty}), \\ [.2cm]
w^{\alpha }(\rR, 0) \geq u_{0}, & \text{for }\left\vert x\right\vert =\rR,
\end{array}%
\right. 
$$
with $\rM_{\varepsilon}>0$. We recall that, usually, the comparison principle is stated under the
assumption that $\alpha \in \rL^{1}(0,\rT)$ for any $\rT\in (0,\rT_{\infty})$ and that in our case we only have  $\alpha \in\rW^{-1,q\prime }(0,\rT)$ for some $q>1$. Nevertheless, the comparison principle still holds since  there exists $A\in\rL^{q}(0,\rT)$, which, as in the proof of
Theorem 1, we can also assume also in $\cC\big ([\rT_{\infty}-\varepsilon, \rT_{\infty})\big )$, for some $\delta >0,$
such that 
\begin{equation*}
\alpha (t)=\frac{d}{dt}A(t).
\end{equation*}%
Thus, to justify the comparison principle we can argue in a similar manner to the case of some stochastic parabolic equations with an additive noise (see, $e.g.$, \cite{DiDi-hetzer}). We make a change of variables which is 
$$
\rU^{\alpha }(t,x)=u^{\alpha }(t,x)-\rA(t).
$$
Then
$$
\left\{ 
\begin{array}{ll}
-\Delta \rU^{\alpha }+\cG(\big (\rU^{\alpha }+\rA(t),\alpha \big )=0 & \text{for }(x,t)\in
\bB_{\rR}\times (0,\rT_{\infty }), \\ [.15cm]
\dfrac{\partial \rU^{\alpha }}{\partial t}+\dfrac{\partial \rU^{\alpha }}{\partial
\bn}=\lambda f_{\rM_{\varepsilon }}(\rU_{\alpha }+A(t)) & \text{for }(x,t)\in
\partial \bB_{\rR}\times (0,\rT_{\infty }), \\ [.2cm]
\rU_{\alpha }(\rR,0)=u_{0}, & \text{for }x\in \partial \bB_{\rR}.%
\end{array}%
\right. 
$$
Analogously, if we define 
$$
\rW^{\alpha }(t)=\rV^{\alpha }(t)-\rA(t)
$$
we get that 
$$
\left\{ 
\begin{array}{ll}
-\Delta \rW^{\alpha }+\cG(\rW^{\alpha }+\rA(t),\alpha \big )\geq 0 & \text{for }(x,t)\in
\bB_{\rR}\times (0,\rT_{\infty }), \\ [.15cm]
\dfrac{\partial \rW^{\alpha }}{\partial t}+\dfrac{\partial \rW^{\alpha }}{\partial
\rn}=\lambda f_{\rM_{\varepsilon }}(\rW_{\alpha }+\rA(t)) & \text{for }(x,t)\in
\partial \bB_{\rR}\times (0,\rT_{\infty}), \\ [.2cm]
\rW^{\alpha }(\rR, 0)\geq u_{0}, & \text{for }x\in \partial \bB_{\rR}.%
\end{array}%
\right. 
$$
Now the comparison principle can be applied (even if the right hand side is a time depending nonlinear term: see, $e.g.$, \cite{DiDi-hetzer}) and we get
that $\rU^{\alpha }(x,t)\leq \rW^{\alpha }(x,t)$ , which implies the desired comparison.
\par
\medskip
\noindent {\sc End of the proof of Theorem 3}. 
For $t\in (0,\rT_{\infty})$ the dynamic boundary condition can be equivalently
expressed as 

\begin{equation*}
\dfrac{\partial u}{\partial t}(\rR,t)=\lambda f_{\rM_{\varepsilon
}}(u(\rR,t))-c(t)+\alpha (t)
\end{equation*}%
where 
\begin{equation*}
c(t)=\dfrac{\partial u}{\partial r}(\rR,t).
\end{equation*}%
Now, we take $\alpha (t)$ as in the proof of Theorem 1. In fact, if  $t\in
(\rT_{\infty}-\varepsilon ,\rT_{\infty}),$ since $\cG(u,\alpha )$ is
truncated, \ we deduce from the identity \eqref{eq:radialmultiplied} and the above comparison Lemma that
\begin{equation*}
0\leq c(t)\leq \rK(\rR)\sqrt{\rV^{\alpha }(t)}
\end{equation*}%
for some $\rK(\rR)>0.$ Then, from the form of the control  $\alpha (t)$ we
arrive to the non-homogeneous neutral equation
$$
\left\{ 
\begin{array}{c}
\dfrac{d}{dt}\left[ \ry(t)-\rB(t)(\ry(t-\tau ))\right] =\lambda
f_{\rM_{\varepsilon }}(\ry(t))-\rB(t)\dfrac{d}{dt}\left[ \ry(t-\tau )\right] -c(t) \\%
[0.2cm]
\ry(\theta )=u^{0}(\theta ),\quad 0\leq \theta \leq \rT_{\infty }-\varepsilon ,%
\end{array}%
\right. 
$$
where now $\ry(t)=u^{0}(\rR,t).$ Then we can apply the generalized {\em Alekseev nonlinear variation of constants formula }and get that its
unique solution $z(t)$ admits the integral representation 
$$
z(t)=\ry^{0}(t)+\rB(t)u^{0}(\rR, t-\tau )-\int_{0}^{t}\big (\rB(s)\dfrac{d}{ds}\big [ \Phi
(t,s,z(s))u^{0}(\rR,s-\tau )\big ]-c(s)\big) ds,
$$
where, again,  $\ry^{0}(t)=\phi (t,0,u_{0})$. \ Then, since the singularity of 
$c(t)$ is weaker than the one of $\rV^{\alpha }(t)$ (which, for this control $%
\alpha (t)$, is of the form $\dfrac{a}{\left\vert \rT_{\infty}-t\right\vert
^{\gamma }}$ with $\gamma \in (0,1)$) we conclude that $z\in
\rL^{1}(0,\rT_{\infty})$). In consequence, since $0\leq u^{0}(r,t)\leq u^{0}(\rR,t)$ we have
that $u\in \rL^{1}(0,\rT_{\infty}:\rL^{\infty}(\bB_{\rR}))$ and the extension is integrable
in the whole domain.  The extension to the rest of the interval $(\rT_{\infty},+\infty)$ is similar and follows as steps 2 and 3 of the proof of Theorem 1.

\begin{rem}\rm \label{Other problems}
The controlled explosions can be also shown for other different partial differential problems leading to a global blow up time ($i.e.$ with the region of blow up given by the entire spatial domain). That was presented in \cite{CDV2008} for the case of the usual semilinear heat equation with Neumman boundary conditions 
$$
\left\{ 
\begin{array}{ll}
\dfrac{\partial u}{\partial t}-\Delta u=\rF(u,\alpha) & \text{for }
(x,t)\in \Omega \times (0,+\infty  ), \\ [.3cm]
\dfrac{\partial u}{\partial \bn}=0, & \text{for }(x,t)\in \partial \Omega \times (0,+\infty   ), \\ [.25cm]
u(0,x)=u_{0}(x), & \text{for }x\in \Omega ,
\end{array}%
\right.   
$$
where $\Omega $ is a regular open bounded set of $\RR^{\rN},~\rN\geq 1$ and $\rF(u,\alpha)$ is associated to a linear delayed term (for the case of separable solutions: Section 4.1 of \cite{CDV2008}) or nonlinear delayed term (Section 4.2 of \cite{CDV2008}). 
A different problem, with a global blow up time was considered in 
\cite{CaDiVe2013}. The formulation was very similar to problem $\rP(\alpha)$ but with a linear elliptic equation ($g=0$). As a matter of fact, this problem can be understood as an special case of the fractional semilinear heat equation:
$$
\quad \left\{ 
\begin{array}{ll}
\dfrac{\partial u}{\partial t}+(-\Delta )^{s}u=\rF(u,\alpha)& 
\text{for }(x,t)\in \Omega \times (0,+\infty  ), \\ [.19cm]
u=0, & \text{for }t\in (0,+\infty  )\text{ and }x\notin \Omega , \\ [.15cm]
u(0,x)=u_{0}(x), & \text{for }x\in \Omega ,%
\end{array}%
\right. 
$$
where $(-\Delta )^{s}u$ represents the fractional Laplacian on $\Omega $,
for $s\in (0,2),$ in the sense of Caffarelli and Silvester \cite{Caff-Silv}. A pioneering 1975 paper dealing with the blow-up question when $\Omega =\RR^{\rN}$ was \cite{Sugitani75}. For the case of $\Omega $ bounded see \cite{Fu-Pucci}. It seems possible to apply the results of Section 2 of this paper to get some controlled explosions but now stated in some different terms. 
We recall that the fractional Laplacian operator when the spatial domain is the whole space $\RR^{\rN}$,  by means of the formula 
$$
(-\Delta )^{s}u(x)=c_{\rN,s}{\rm P.V.}\int_{\RR^{\rN}}\frac{u(x)-u(y)}{
|x-y|^{\rN+2s}}dy,  \label{frac.lap}
$$
with parameter $s\in (0,1)$ and a precise constant $c_{\rN,s}>0$ that we do
not need to make explicit for our purposes. The operator can also be defined
via the Fourier transform on $\RR^{\rN}$. With an appropriate value
of the constant $c_{n,s}$, the limit $s\nearrow 1$ produces the classical
Laplace operator $-\Delta $, while the limit $s\rightarrow 0$ is the
identity operator. An equivalent definition of this fractional Laplacian
uses the so-called extension method, that was well known for $s=\frac{1}{2}$
and has been extended to all $s\in (0,1)$ by \cite
{Caff-Silv}. We recall that the existence,
comparison principle (implying the uniqueness) for local in time solutions
for this type of problems are consequence of well known results (see, e.g., \cite{Vazquez-Cime, Sugitani75, Fu-Pucci} and their many references). 
It can be shown (see, $e.g.$, \cite{Dgc-Vazquez}) that a good space to solve this problem is $
\rX=\rL^{1}(\Omega :\delta )=\big \{ w\in\rL_{loc}^{1}(\Omega ):\delta w\in
\rL^{1}(\Omega )\big\} $, with $\delta (x)=\hbox{dist}(x,\partial \Omega)$. The
adaptation, to this problem, of many of the results of \cite{Brezis-Cacenave} is
automatic and then if $f$ satisfies (\ref{Hypo superlinear}), $f$ is convex
and $f(0)>0,$ then there exists a $\lambda ^{\ast }>0$ such that the local
very weak solutions blows up in finite time if and only if $\lambda >\lambda
^{\ast }$. Then, it is possible to show that if
$u_{0}\in \rL^{\infty }(\Omega ),$ $u_{0}\geq 0$, and $\lambda >0$ is such that the local very weak solution $u^{0}(x,t)$ of problem with no control, blows up in a finite time $\rT_{\infty }$. The proof of the control of the the trajectory $u^{0}(\cdot,t)$ in a sense quite similar to the one of Definition 1 will be given in a separated work. $\fin$
\end{rem}
\bigskip 

\begin{center}
{\bf \sffamily Acknowledgements}	
\end{center}
The research of A.C. Casal, G. D\'{\i}az and J.I. D\'{\i}az was partially supported
by the project ref. PID2020-112517GB-I00 of the Ministerio de Ciencia e Innovaci\'{o}n
(Spain).

{\footnotesize
\begin{tabular}{llll}
Alfonso Carlos Casal & Gregorio D\'{\i}az & Jesús Ildefonso  D\'{\i}az & Jos\'{e} Manuel Vegas\\
	& & Instituto Matem\'{a}tico & \\
	& & Interdisciplinar (IMI) & \\
	Dpto. de Matem\'{a}tica Aplicada & Dpto. An\'{a}lisis Matemático & Dpto. An\'{a}lisis Matemático & \\
	& y Matem\'atica Aplicada & y Matem\'atica Aplicada & \\
	ETS de Arquitectura & U. Complutense de Madrid & U. Complutense de Madrid & Cunef \\
		28040 Madrid, Spain& 28040 Madrid, Spain & 28040 Madrid, Spain &	28040 Madrid, Spain.\\
	{\tt  alfonso.casal@upm.es}& {\tt gdiaz@ucm.es} & {\tt jidiaz@ucm.es} & {\tt jm.vegas@cunef.edu} 
\end{tabular}
}


\begin{thebibliography}{99}                                                                                               %
\bibitem {AlDiRe}S. Alarc\'{o}n, G. D\'{\i}az, G. and J.M. Rey: Large solutions
of elliptic semilinear equations in the borderline case. An exhaustive and
intrinsic point of view, {\em Journal of Mathematical Analysis and
Applications}, {\bf 431} (2015), 365-405.

\bibitem {Alekseev}V.M. Alekseev: An estimate for the perturbations of the
solutions of ordinary differential equations (Russian), {\em Vestnik Moskov
Univ. Ser. I Mat. Meh.}, {\bf 2} (1961), 28--36.

 \bibitem {Amman-Fila}H. Amann and M. Fila: A Fujita-type theorem for the
Laplace equation with a dynamical boundary condition, {\em Acta Math. Univ.
Comenian.}, {\bf 66} {\bf 2} (1997), 321--328.

\bibitem{Amman-Quittner}H.Amann and P. Quittner: Optimal control problems governed by semilinear parabolic equations
with low regularity data, {\em Adv. Differ. Equ.} {\bf 11}, (2006) 1-33.

\bibitem {Arrieta}J.M. Arrieta: On boundedness of solutions of
reaction-diffusion equations with nonlinear boundary equations, {\em Proc.
Amer. Math. Soc.}, {\bf 136} {\bf 1} (2007), 151-160.



\bibitem {Balch}K. Le Balc'h: Global null-controllability and nonnegative-controllability of slightly superlinear
heat equations, {\em J. Math. Pures Appl.} {\bf 9} 135 (2020), 103-139.

\bibitem {Bandle}C. Bandle: Asymptotic behavior of large solutions of
elliptic equations, {\em Annals of University of Craiova, Math. Comp. Sci.
Ser.}, {\bf 32} (2005), 1-8.

\bibitem {BaDiDi}C. Bandle, G. D\'{\i}az et J.I. D\'{\i}az: Solutions d'
Équations de r\'{e}action-diffusion non lin\'{e}aires explosant
au bord parabolique, {\em C.R.Acad Sci Paris}, {\bf 318} (1994), 455-460.


\bibitem {Bandle-Reichel}C. Bandle, J. von Below  and W. Reichel: Parabolic
problems with dynamical boundary conditions: eigenvalue expansions and blow
up, {\em Atti della Accademia Nazionale dei Lincei}, Classe di Scienze
Fisiche, Matematiche e Naturali, Rendiconti Lincei Matematica E Applicazioni
2006, 35-67.

\bibitem {Baras-Cohen}P. Baras  and L. Cohen: Complete Blow-Up after
$\rT_{\max}$ for the Solution of a Semilinear Heat {\em Equation, J.
Funct. Anal.} {\bf 71} (1987), 142-174.


\bibitem {BeDiVr}I. Bejenaru, J.I. D\'{\i}az and I.I. Vrabie: An abstract
approximate controllability result and applications to elliptic and parabolic
systems with dynamical boundary conditions. {\em Electr. J. Diff. Eqns.},
{\bf 50} (2001), 1-19.


\bibitem {BreOMM}H. Brezis: {\em Operateurs maximaux monotones et
semi-groupes de contractions dans les espaces de Hilbert}, North-Holland
Mathematical Studies, Amsterdam, 1973.


\bibitem {Brezis-Cacenave}H. Brezis,Th. Cazenave, Y. Martel and A.
Ramiandrisoa: Blow up for $u_{t}-\Delta u=g(u)$ revisited. {\em Adv.
Differ. Equat.}, {\bf 1} (1996), 73--90.

\bibitem {Cacenave-Matel-Zhao}Th. Cazenave, Y. Martel and L. Zhao:
Solutions blowing up on any given compact set for the energy subcritical wave
equation. \textit{J. Differential Equations} {\bf 268} (2020), no. 2, 680--706.

\bibitem{Caff-Silv} L. Caffarelli and L. Silvestre: An extension problem
related to the fractional laplacian, {\em Commun. Partial Differential Equations} 
\textbf{32} (8) (2007) 1245--1260.

\bibitem{CDV2008} A.C. Casal, J.I. D\'{\i}az, and J.M. Vegas: Blow-up in
some ordinary and partial differential equations with time-delay. {\em
Dynam. Systems Appl.} {\bf 18} (1) (2009), 29--46.

\bibitem {CaDiVe2013}A.C. Casal, J.I. D\'{\i}az and J.M. Vegas: Controlled
explosions of blowing-up trajectories in semilinear problems and a nonlinear
variation of constant formula, {\em XXIII Congreso de Ecuaciones
Diferenciales y Aplicaciones, XIII Congreso de Matem\'{a}tica Aplicada},
Castell\'{o}n, 9-13 septiembre 2013. e-{\em Proccedings. }



\bibitem {CaDiVe}A.C. Casal, J.I. D\'{\i}az and J.M. Vegas: Complete
recuperation after the blow up time for semilinear problems, {\em AIMS
Procceding 2015}, 223-229.



\bibitem {Coron-Trelat}J.M. Coron and E. Tr\'{e}lat: Global steady-state
controllability of 1-D semilinear heat equations, {\em SIAM J. Control and
Optimization}, {\bf 43} (2) (2004), 549-569.


\bibitem {DiLe}G. D\'{\i}az and R. Letelier: Explosive solutions of
quasilinear elliptic equations: existence and uniqueness, {\em Nonlinear
Analysis}, {\bf 20} (1993), 97-125.


\bibitem{DDO} G. D\'{\i}az, J.I. D\'{\i}az and J. Otero:  Construction of the maximal solution of Backus’ problem in geodesy and geomagnetism. Studia Geophysica et Geodaetica, 55(3) (2011), 415-440.

\bibitem {DiDi-hetzer}G. D\'{\i}az  and J.I. D\'{\i}az: Stochastic energy balance climate models with Legendre weighted diffusion and a cylindrical Wiener process forcing, {\em Discrete and Continuous Dynamical Systems  Series S}. doi: 10.3934/dcdss.2021165.


\bibitem{Dgc-Vazquez}J.I. D\'{\i}az, D. Gómez-Castro, and J. L. Vázquez. The fractional Schrödinger equation with general nonnegative potentials. The weighted space approach. Nonlinear Analysis, Nonlinear Analysis, 177 (2018) 325-360

\bibitem {Diaz-Lions}J.I. D\'{\i}az and J.L. Lions: Sur la
contr\^{o}labilit\'{e} de probl\`{e}mes paraboliques avec phenomenes
d'explosion, {\em C. R. Acad. Scie. de Paris}. t. 327, S\'{e}rie I (1998), 173-177.

\bibitem {Diaz-Lions 2}J.I. D\'{\i}az and J.L. Lions: On the approximate
controllability for some explosive parabolic problems. In: Hoffmann, K.-H., et
al. (eds.) {\em Optimal Control of Partial Differential Equations} (Chemnitz,
1998), Internat. Ser. Numer. Math., vol. 133 (1993), Birkh\"{a}user, Basel, 115--132.




\bibitem {FdezCara-Zuazua}E. Fern\'{a}ndez-Cara and E. Zuazua: Null and
approximate controllability for weakly blowing up semilinear heat equations.
{\em Ann. Inst. Henri Poincar\'{e} Anal. Non Lin\`{e}aire}, {\bf 17} (2000),
583--616.

\bibitem {Fila-Filo}M. Fila and J. Filo: Blow-up on the boundary: a survey.
In {\em Singularities and differential equations}. Banach Center
Publications, volume 33 ()1996). Institute of Mathematics, Polish Academy of Sciences,
Warszawa, 67-77.

\bibitem {Fila-Quitner}M. Fila and P. Quittner: Large Time Behavior of
Solutions of a Semilinear Parabolic Equation with a Nonlinear Dynamical
Boundary Condition, in Topics in Nonlinear Analysis: The Herbert Amann
Anniversary Volume, (Joachim Escher and Gieri Simonett eds.), Progress in
Nonlinear Differential Equations and Their Applications, Vol. 35 (1999), Birkhauser, 251-272.



\bibitem{Fu-Pucci}Y. Fu and P. Pucci: On solutions of space-fractional diffusion equations by means of potential wells. Electronic Journal of Qualitative Theory of Differential Equations, 2016 (70), 1-17.



\bibitem {Galaktionov-Vazquez Libro}V.A. Galaktionov and J.L. V\'{a}zquez: A stability technique for evolution partial differential equations: a dynamical systems approach.  Progress in Nonlinear Differential Equations and Their Applications Vol. 56. Birkhauser, Boston, 2003.



\bibitem {Hu}B. Hu: {\em Blow-up Theories for Semilinear Parabolic
Equations}, Lecture Notes in Mathematics 2018, Springer-Verlag, Berlin, 2011.


\bibitem {Kirane}M. Kirane: Blow-up for some equations with semilinear dynamical boundary conditions of parabolic and hyperbolic type. {\em Hokkaido Math. J.},
{\bf 21} 2 (1992), 222-229.

\bibitem {KiraneNabanaPokhozhaev}M. Kirane, E. Nabana and S.I. Pokhozhaev: The Absence of Solutions of Elliptic Systems with Dynamic Boundary
Conditions, {\em Differential Equations}, {\bf 38} 6 (2002), 808--815.


\bibitem {LaksLeela} V. Laksmikantham and S. Leela: {\em Differential and
Integral Inequalities, Theory and Applications}, Vols. I and II, Academic
Press, New York, 1969.

\bibitem {Levine-Payne}H. Levine and L. Payne: Nonexistence theorems for the
heat equations with nonlinear boundary conditions and for the porous medium
equation backward in time, {\em J. Differential Equations}, {\bf 16} (2) (1974),
319-334.



\bibitem {Lions}J.L. Lions: {\em Contr\^{o}le des syst\^{e}mes
distribu\'{e}s singuliers}, Gauthier-Villars, Bordas, Paris. 1983.

\bibitem {Lopez-Wolansky}J. L\'{o}pez G\'{o}mez, V. M\'{a}rquez and N.I. Wolanski: Blow up results and localization of blow up points for the heat equation
with a nonlinear boundary condition, {\em J. Diff. Eq.} {\bf 92} (2) (1991), 384-401.


\bibitem {Merle}F. Merle: Solution of a nonlinear heat equation with
arbitrarily given blow-up points. {\em Comm. Pure Appl. Math}. {\bf 45} no. 3
(1992), 263--300.


\bibitem {Porreta-Zuazua}A. Porretta and E. Zuazua: Null controllability of viscous Hamilton-Jacobi equations, {\em Ann. Inst. H. Poincaré Anal. Non Linéaire} {\bf 29} (2012), 301-333.

\bibitem {Quitner-Souplet}P. Quittner and P. Souplet: {\em Superlinear
Parabolic Problems}, Birkh\"{a}user, Berlin, 2007.




\bibitem{Sugitani75} S. Sugitani: On nonexistence of global solutions for some nonlinear integral equations. Osaka Journal of Mathematics, 12(1) (1975), 45-51.



\bibitem{Vazquez-Cime} J. L. V\'{a}zquez. The mathematical theories of diffusion. Nonlinear and fractional diffusion, Lecture Notes in Mathematics,
2186. Fondazione CIME/CIME Foundation Subseries. Springer, Cham; Fondazione
C.I.M.E., Florence, 2017.

\bibitem {Vazquez-Vitillaro}J.L. V\'{a}zquez and E. Vitillaro, E.: On the
Laplace equation with dynamical boundary conditions of reactive
diffusive type, {\em J. Math. Anal. Appl.}, {\bf 354} 2 (2009), 674-688.



\bibitem {Ya} N. Yamazaki, A class of nonlinear evolution equations governed by time-dependent operators of subdifferential type, Hokkaido University Preprint Series in Mathematics 696 (2005), 1-16.


\end{thebibliography}
\end{document}